\documentclass[A4paper,11pt]{article}
\usepackage{amssymb}
\usepackage{amsmath}
\usepackage{theorem}
\usepackage{graphicx}
\usepackage{mathrsfs}
\usepackage[all]{xy}

\theoremstyle{break}
\theorembodyfont{\normalfont}
\newtheorem{Def}{Definition}[section]
\newtheorem{Prop}[Def]{Proposition}
\newtheorem{Cor}[Def]{Corollary}
\newtheorem{Lem}[Def]{Lemma}

\newtheorem{Thm}[Def]{Theorem}
\newtheorem{Rmk}[Def]{Remark}
\newtheorem{Hyp}[Def]{Hypothesis}
\newtheorem{Pf}{Proof}

\newcommand{\bfvec}[1]{\mbox{\boldmath $#1$}} 

\newcommand{\fd}[1]{\mathcal{O}(\mathcal{M}_{A/#1 I, F})}
\newcommand{\sie}[1]{\mathrm{g}_{#1}}
\newcommand{\Kinfty}{K_\infty^{(p)}}
\newcommand{\Gal}[2]{\mathrm{Gal}\, (#1/#2)}

\title{Constant Terms of Coleman Power Series and Euler Systems in Function Fields}
\author{Toshiya SEIRIKI}
\date{\today}

\begin{document}
\maketitle

\begin{abstract}
In this paper, we calculate the constant terms of Coleman power series, and give a proof of Iwasawa main conjecture over global function fields of characteristic $p > 0$ using Euler systems. The calculation of the constant terms of Coleman power series is necessary to prove a relation of Kolyvagin's derivative classes induced by our Euler systems. This paper is a revised version of the author's master's thesis.
\end{abstract}

\section{Introduction} 
\label{intro}

In 1978 Robert Coleman established in his original paper \cite{Co} the theory of norm compatible systems for the tower of cyclotomic extensions over complete discrete valuation fields whose residue fields are finite. Let $K$ be a complete discrete valuation field whose residue field has $q$ elements and $\mathcal{O}_K$ the integer ring of $K$. Fix a uniformizer $\pi$ of $\mathcal{O}_K$. Let $\mathcal{F}$ denote the Lubin-Tate $\mathcal{O}_K$-module corresponding to $\pi$, $[\pi]_F$ the endomorphism of $\mathcal{F}$ given by $\pi$, which satisfies 
\[
  [\pi]_F \equiv \pi X \mod \mathrm{deg}\, 2 \quad , \quad  [\pi]_F \equiv X^q \mod \pi , \quad N_F(X) = X
\]
where $N_F$ is the Coleman's norm operator. Let $\mathcal{F}(n)$ be the set of roots of $\pi^{n+1}$ in $\mathcal{F}$, $H$ denote a maximal unramified extension of $K$ and $\mathrm{Frob}$ the arithmetic Frobenius of $\Gal{H}{K}$. We define the tower of extensions: 
\[
	H_n := H(\mathcal{F}(n))\ (n=0,1,2,\ldots).
\]
Let $\xi = (\xi_n)_n$ be a generator of $\varprojlim_n \mathcal{F}(n) $ as an $\mathcal{O}_K$-module. Then for any $m \geq n \geq 0$, we have
\[
	[\pi^{m-n}]_F(\xi_m) = \xi_n.
\]
Let $\mathcal{O}_H[[T]]$ denotes the ring of formal power series whose coefficients are in $\mathcal{O}_H$. The group $\Gal{H}{K}$ acts on $\mathcal{O}_H[[T]]$ by the action on the coefficients. Then for any $\alpha = (\alpha_n)_n \in \varprojlim \mathcal{O}_{H_n}^\times$, where the inverse limit is taken with respect to norm maps, there exists a unique element $ \mathrm{Col}_\alpha(T)$ in $\mathcal{O}_H[[T]]^\times$ such that
	\[
		(\mathrm{Frob}^{-n}\mathrm{Col}_\alpha)(\xi_n) = \alpha_n.
	\]
We call $\mathrm{Col}_\alpha(T)$ the Coleman power series for $\alpha$. When an element $\alpha = (\alpha_n)_n$ is given, it is difficult to compute the Coleman power series for $\alpha$. However it is known that, for example, we can prove the analogue of Wiles' explicit reciprocity law for rank 1 Drinfeld modules using Coleman power series and ``Dlog'' (see \cite{BL}). When we use the operator Dlog, the information of constant term disappears.

Our first result of this paper is to calculate the constant terms of Coleman power series.  We give another way to analysis the power series. In section \ref{CTCPS} we will prove the following theorem.

\begin{Thm}\label{maintheorem1}
	Let $K, H, H_n$ be as above. Let $u = (u_n)_n$ be an element in $\varprojlim_n \mathcal{O}_{H_n}^\times$, $\mathrm{Col}_u(T)$ the Coleman power series for $u$. We assume that $u_H = 1$. Put $H_\infty := \underset{n}\cup H_n$. Then for any character $\chi : \Gal{H_\infty}{H} \rightarrow \mathbb{Q}/\mathbb{Z}$ whose order is finite, the following equality holds:
	\[
		\mathrm{rec}_{H}^{-1}((u_{H_\chi}, \cdot)_{H_\chi/H}) = \mathrm{Col}_u(0).
	\]
Here $H_\chi := (H_\infty)^{\mathrm{Ker}\chi}$ is an intermediate field of $H_\infty / H$ corresponding with $\mathrm{Ker}\, \chi$, $(u_{H_\chi}, \cdot)_{L/K}$ is a symbol which is introduced in section \ref{symb} and $\mathrm{rec}_{H}$ is the reciprocity map.
\end{Thm}

Our second result is an analogue of Iwasawa main conjecture over global function fields of characteristic $p>0$.
Let $v$ be a prime of a global function field $F$, $K_n := \mathcal{O} ( \mathcal{M}_{A/v^nI, F} ), n \geq 0$ extensions of $F$ (see section \ref{seceu} for the definition of the notation) and $\Delta$ the Galois group of $\mathcal{O} ( \mathcal{M}_{A/vI, F} )/F$.
These extensions $\mathcal{O} ( \mathcal{M}_{A/v^nI, F} )$ are the analogues of tower of cyclic extensions of $\mathbb{Q}$. 
For any $p$-adic character $\rho$ in $\Delta$, put 
\[
	e(\rho) := \frac{1}{\sharp\Delta}\displaystyle\sum_{\delta\in\Delta}\, \rho (\delta)^{-1}\delta.
\]
Fix $K_\infty^{(p)}$ an intermediate field of $K_\infty / K_0$ such that $\Gamma := \Gal{K_\infty^{(p)}}{K_0} \cong \mathbb{Z}_p$.
Put $K_n^{(p)} := (K_\infty^{(p)})^{\Gamma_n}$ where $\Gamma_n := \Gamma^{p^n}$.
Let $C_n$ denote the p-part of ideal class group of $\mathcal{O} ( \mathcal{M}_{A/v^nI, F} )$, $E_n$ the units group, $\mathcal{E}_n \subset E_n$ the group generated by Siegel units, $\overline{E_n}$ and $\overline{\mathcal{E}_n}$ the closure of $E_n \cap U_n$ and $\mathcal{E}_n$ and $X_\infty$ the inverse limit $\varprojlim X_n$ for $X = C, E, \overline{\mathcal{E}}$. Under some technical assumptions (see hypothesis \ref{hyp0}) we can prove the following theorem:

\begin{Thm}[Iwasawa main conjecture]
	Suppose that p $\nmid \#\Delta$. For any irreducible character $\rho$ in $\Delta$,
	\begin{align*}
		{\rm char}\, e(\rho)(C_\infty) = {\rm char}\left(e(\rho)E_\infty)/e(\rho)\overline{\mathcal{E}_\infty}\right)
	\end{align*}
	where ${\rm char}(X)$ is the characteristic ideal of X. 
\end{Thm}

We describe a structure of ideal class groups of cyclic extensions using an analogue of Rubin's method over a cyclotomic number fields. His idea in the classical case is using Euler systems constructed by cyclotomic units and bounding the number of ideal class groups (see \cite{R1}). In our settings, our Euler systems are constructed by Siegel units (see section \ref{defsie} and \ref{seceu}).
 and $M$ a power of $p$. We consider Euler systems in $\mathcal{O} ( \mathcal{M}_{A/v^nI, F} )$ and their Kolyvagin's derivative classes $\{ \kappa_{(\bfvec{v}_r ,  \psi_{\bfvec{v}_r})} \}_r$ (see section \ref{kap}). Then the following equalities named ``finite-singular comparison equalities'' is the key point of proof of Iwasawa main conjecture.
\begin{Thm}[finite-singular comparison equalities]\label{maintheorem2}
	\begin{align*}
		\quad [\kappa_{(\bfvec{v}_r ,  \psi_{\bfvec{v}_r})}]_v
		 = \underset{\lambda \mid v}\sum \psi_v 
			(\kappa_{(\bfvec{v}_{r \backslash i},\psi_{\bfvec{v}_{r \backslash i}})})\lambda
	\end{align*}
	where $[\cdot]_v$ is a projection of a $v$-part of a principal ideal onto $\mathscr{I}_v/M\mathscr{I}_v$, $\mathscr{I}_v$ is the $v$-part of group of fractional ideals of $\mathcal{O} ( \mathcal{M}_{A/I, F} )$.
\end{Thm}
When $F$ is a number field, we get this equality due to Kummer theory since the characteristic of $F$ is zero. However in our case we cannot apply Kummer theory, so we have to find out another way to prove this equality. In this paper we find out a method with which we can apply the theory of Euler systems to such a case by using the theory of the Coleman power series. We prove that the constant term of Coleman power series for a Siegel unit is equal to another Siegel unit (see Proposition \ref{20}).

We remark that if $M$ is prime to $p$ then Theorem \ref{maintheorem2} was proved by Hassan Oukhaba and St$\acute{{\rm e}}$phane Vigui$\acute{{\rm e}}$ by a similar method. However their method cannot be applied to the case $M$ is a power of $p$. 

The author would like to thank Seidai Yasuda and Satoshi Kondo for his helpful comments.


\section{Calculate the constant terms of Coleman power series}\label{constterm}
In this section we prove Theorem \ref{maintheorem1}. First we introduce a new symbol $(\cdot, \cdot)_{L/K}$. The constant terms of Coleman power series is characterized by this symbol. 

\subsection{The symbol $(\cdot ,\cdot )_{L/K}$}\label{symb}
Let $K$ be a complete discrete valuation field whose residue field is finite, $L$ a finite abelian extension of $K$, $d$ the degree of the extension, $G$ the Galois group $\Gal{L}{K}$, $N_{L/K}$ the norm map from $L$ to $K$ and $\mathrm{val}_K$ a valuation of $L$
which is normalized at $K$, i.e., $\mathrm{val}_K(K^\times) = \mathbb{Z}$. Put $U_1 := \mathrm{Ker}\ N_{L/K}$.  
Let $\widehat{G}$ denote the group of characters of $G$ with values in $\mathbb{Q}/\mathbb{Z}$. For any $\chi \in \widehat{G}$, let $K_\chi$ be the intermediate field of $L/K$ which corresponds to $\mathrm{Ker}\, \chi$, i.e. $\Gal{L}{K_\chi} = \mathrm{Ker}\, \chi$. It is clear that $K_\chi$ is a cyclic extension of $K$. Write $d_\chi$ for the degree of the extension and $G_\chi$ its Galois group. For $u \in U_1$, let $u_{K_\chi}$ denote the image of $u$ under the norm map $N_{L/K_\chi}$.  
Let us choose $\sigma \in G$ which satisfies $\chi(\sigma) = 1/m$ where $m$ is the order of $\chi$. We can show that $N_{K_\chi/K}(u_{K_\chi}) = 1$ since $u \in U_1$. Therefore using Hilbert 90, there exists $b_\chi \in K_\chi^\times$ such that
\[
	u_{K_\chi} = \frac{\sigma(b_\chi)}{b_\chi}.
\]
It is unique up to $K^\times$. 

\begin{Def}
	We define the map $(\cdot , \cdot)_{L/K} : U_1 \times \widehat{G} \rightarrow \mathbb{Q}/\mathbb{Z}$ as 
	\[
		(u , \chi)_{L/K} = \mathrm{val}_K(b_\chi) \ \mathrm{mod}\ \mathbb{Z}
	\]
	for any 	$u \in U_1 , \chi \in \widehat{G}$. 
\end{Def}

This is well-defined since $\mathrm{val}_K$ is normalized, and it is clear that
\[
	(uu' , \chi)_{L/K} = (u , \chi)_{L/K} + (u' , \chi)_{L/K}.
\]
However it is not easy to show additivity of the other side.

\begin{Prop}\label{bilinear}
	Assume that $L/K$ is a totally ramified abelian extension. Then for any $u \in U_1$ and elements $\chi,\chi' \in \widehat{G}$, the following equality holds:
	\[
		(u , \chi \chi')_{L/K} = (u , \chi)_{L/K} + (u , \chi')_{L/K}. 
	\]
\end{Prop}

\begin{Pf}
	Fix $\pi_L \in L$ a uniformizer. Then for any $u \in U_1$, there exist a finite unramified extension $K'/K$, units $\beta_1,\ldots,\beta_r \in \mathcal{O}_{K'}^\times$ and $\tau , \sigma_1,\ldots,\sigma_r \in \Gal{LK'}{K'}$ such that
	\[
		u = \displaystyle\frac{\tau\left( \pi_L \right)}{\pi_L}\prod_{i=1}^r \frac{\sigma_i\left(\beta_i\right)}{\beta_i}
	\]
	since 
	\begin{align*}
		\underset{K'/K:unram}\varinjlim \widehat{H}^{-1}(G,(LK')^\times) &\cong H_2(G, \mathbb{Z}) \\
		&\cong \wedge^2 G = 0.
	\end{align*}
	Therefore we have
	\begin{align*}
		(u , \chi)_{L/K} &= \left( \frac{\tau(\pi_L)}{\pi_L} ,\chi \right)_{LK'/K'} +
			\displaystyle\sum_{i=1}^r \left( \frac{\sigma(\beta_i)}{\beta_i} , \chi \right)_{LK'/K'} \\
		&=  \left( \frac{\tau(\pi_L)}{\pi_L} , \chi \right)_{LK'/K'}.
	\end{align*}
	By the definition, the last term is equal to $(\chi)(\tau)$ for any $\chi \in \widehat{G}$. It implies
	\[
		(u , \chi \chi')_{L/K} = (u , \chi)_{L/K} + (u , \chi')_{L/K}. 
	\]
\end{Pf}

\begin{Prop}\label{symbol0}
	Fix $u \in U_1$ and $\chi \in \widehat{G}$. Assume that $L$ has no nontrivial $d_\chi$-th roots of unity, and that
	there exist elements $a \in \mathcal{O}_K^\times , b \in \mathcal{O}_L^\times$ such that $u = ab^{d_\chi}$.
	Then the following equality holds:
	\[
		(u , \chi)_{L/K} = 0. 
	\]
\end{Prop}

\begin{Pf}
	It is enough to give a proof in case $L = K_\chi$. Let $\sigma$ be a generator of $G_{\chi}$.  
	If there exists $u' \in \mathcal{O}_{K_\chi}^\times$ such that $u = \sigma(u')/u'$, then we have $(u , \chi)_{L/K} = 0$
	via the definition. Therefore we have to prove the existence of $u'$. Put 
	\[
		u' := \displaystyle \prod_{i=0}^{d_\chi -1} \sigma^i(b^i)
	\]
	and then we have
	\begin{align*}
		(\sigma(u')/u')^{d_\chi} &= \left( \sigma \displaystyle \prod_{i=0}^{d_\chi-1} \sigma^i(b^{d_\chi i})\right) / 
										\left( \displaystyle \prod_{i=0}^{d_\chi-1} \sigma^i(b^{d_\chi i})\right) \\
							&= \left( \sigma \displaystyle \prod_{i=0}^{d_\chi-1} \sigma^i(a^i b^{d_\chi i})\right) / 
										\left( \displaystyle \prod_{i=0}^{d_\chi-1} \sigma^i(a^i b^{d_\chi i})\right) \\
							&= \left( \sigma \displaystyle \prod_{i=0}^{d_\chi-1} \sigma^i(u^i)\right) / 
										\left( \displaystyle \prod_{i=0}^{d_\chi-1} \sigma^i(u^i)\right) \\
							&= \left( \displaystyle \prod_{i=1}^{d_\chi} \sigma^i(u^{i-1})\right) / 
										\left( \displaystyle \prod_{i=0}^{d_\chi-1} \sigma^i(u^{i-1})\right) \\
							&= \left( \sigma \displaystyle \prod_{i=1}^{d_\chi} \sigma^i(u^i)\right) / 
										\left( \displaystyle \prod_{i=0}^{d_\chi-1} \sigma^i(u^{i-1})\right)
										\left( \displaystyle \prod_{i=0}^{d_\chi-1} \sigma^i(u)\right) \\
							&= \left( \sigma^{d_\chi}(u^{d_\chi}) \right) / 
										\left( \displaystyle \prod_{i=0}^{d_\chi-1} \sigma^i(u)\right) = u^{d_\chi}/1
	\end{align*}
	By the assumption, we can erase $d_\chi$-th power. 
\end{Pf}

\begin{Lem}\label{13}
	Assume that $L$ is a finite totally ramified abelian extension of $K$. Let $N_{L/K}$ denote the norm map.
	For any $u \in \mathcal{O}_K^\times$,
	there exist a finite unramified extension $K'$ and a unit $u' \in \mathcal{O}_{L'}^\times$ where $L' = K'L$ such that
	\[
	u = N_{L'/K'}\left(u'\right)
	\]
	holds.
\end{Lem}

\begin{Pf}
	Since $N_{L/K}\left(\mathcal{O}_L^\times\right)$ is a subgroup of finite index of $\mathcal{O}_K^\times$, 
	there exists an integer
	$n \geq 1$ such that $u^n \in N_{L/K}\left(\mathcal{O}_L^\times\right)$. Let $K'$ be the unramified extension of degree $n$
	of $K$. Then the following homomorphism
	\[
		\mathcal{O}_{K'}^\times / N_{L'/K'}\left(\mathcal{O}_{L'}^\times\right)
			\rightarrow \mathcal{O}_K^\times / N_{L/K}\left(\mathcal{O}_L^\times\right)
	\]
	induced by the norm map $N_{K'/K}$ is an isomorphism. Since we regard $u$ as an element in $\mathcal{O}_{K'}$, we have
	\[
		N_{K'/K}\left(u\right) = u^n \in N_{L/K}\left(\mathcal{O}_L^\times\right).
	\]
	Therefore we have $u \in N_{L'/K'}\left(\mathcal{O}_{L'}^\times\right)$. It implies that there exists a unit
	$u' \in \mathcal{O}_{L'}^\times$ such that $u = N_{L'/K'}\left(u'\right)$ holds.
\end{Pf}

\begin{Lem}\label{14}
	Assume that $L$ is a finite totally ramified abelian extension of $K$ and 
      $u \in \mathcal{O}_L^\times$ satisfies $N_{L/K}\left(u\right) = 1$ and $\left(u,\chi\right)_{L/K} = 0$
	for any character $\chi : \Gal{L}{K} \rightarrow \mathbb{Q}/\mathbb{Z}$. Then there exist a finite 
	unramified extension $K'$, an integer $r \geq 0$, units $\beta_1,\ldots,\beta_r \in \mathcal{O}_{L'}^\times$ and elements
	$\sigma_1,\ldots,\sigma_r \in \Gal{L'}{K'}$ such that the following equality holds:
	\[
		u = \displaystyle\prod_{i=1}^r \frac{\sigma_i\left(\beta_i\right)}{\beta_i}
	\]
       where $L' = K'L$.
\end{Lem}

\begin{Pf}
	We prove it by induction on the number of generators of $\Gal{L}{K}$. When $\Gal{L}{K} = 1$, the claim is clear.
	When $\Gal{L}{K} \neq 1$, there exists an intermediate extension $M$ of $L/K$ such that the number of generators 
	$\Gal{M}{K}$ is less than that of $\Gal{L}{K}$ and $L/M$ is a cyclic extension. Thus there exists an intermediate extension 
	$N$ of $L/K$ such that the composition
	\[
		\Gal{L}{M} \hookrightarrow \Gal{L}{K} \twoheadrightarrow \Gal{N}{K}
	\]
	is a bijection. Let us apply the inductive hypothesis to $N_{L/M}\left(u\right) \in \mathcal{O}_M^\times$. There exist a finite 
	unramified extension $K''$, an integer $r' \geq 0$, units $\beta'_1,\ldots,\beta'_{r'} \in \mathcal{O}_{M''}^\times$ and elements
	$\sigma'_1,\ldots,\sigma'_{r'} \in \Gal{M''}{K''}$ such that 
	\begin{align}\label{01}
		N_{L/M}\left(u\right) = \displaystyle\prod_{i=1}^{r'} \frac{\sigma'_i\left(\beta'_i\right)}{\beta'_i}
	\end{align}
	where $M'' = K''M$. Let us apply Lemma \ref{13} to $\beta'_1,\ldots,\beta'_{r'}$. Then there exist a finite unramified extension
	$K'$ of $K''$ and units $\beta_1,\ldots,\beta_{r'} \in \mathcal{O}_{L'}^\times$ such that
	\[
		\beta'_i = N_{L'/M'}\left(\beta_i\right)
	 \]
	for any $i=1,\ldots,r'$ where $L' = K'L$ and $M' = K'M$. Let $\sigma_1,\ldots,\sigma_{r'} \in \Gal{L'}{K'}$ denote liftings
	of $\sigma'_1,\ldots,\sigma'_{r'} \in \Gal{M''}{K''} \cong \Gal{M'}{K'}$. Put
	\[
		u' := u\cdot\left(\displaystyle\prod_{i=1}^{r'} \frac{\sigma_i\left(\beta_i\right)}{\beta_i}\right)^{-1}.
	\]
	Then we have $N_{L'/M'}\left(u'\right) = 1$, by the equation (\ref{01}). Let us choose a 
	generator $\tau \in \Gal{L'}{M'}$, since by our assumption that $\Gal{L'}{M'} \cong \Gal{L}{M}$ is cyclic. By Hilbert 90 
	there exists an element $b \in L'^\times$ such that $u' = \tau\left(b\right)/b$. Moreover we have 
	$\left(u',\chi\right)_{L'/K'} = 0$ for any character 
	$\chi : \Gal{L'}{K'} \rightarrow \mathbb{Q}/\mathbb{Z}$. Put $m = [L:M]$ and $N' = K'N$. Let us choose $\chi$ so that
	$\chi\left(\tau\right) = 1/m$ and $\chi(\sigma) = 1$ for all $\sigma \in \Gal{L'}{N'}$. This gives us a character on 
	$\Gal{L'}{K'}$ since
	$M'N' = L'$ and $M' \cap N' = K'$. From $(u' , \chi)_{L'/K'} = 0$, we deduce there exists an element $a \in M'^\times$
	such that $ab \in \mathcal{O}_{L'}^\times$. Put $r = r' +1$ and $\beta_r = ab$. Since $a \in M'^\times$ we have 
	$u' = \tau\left(\beta_r\right)/\beta_r$. Therefore, put $\sigma_r = \tau$ and we have 
	\[
		u = \displaystyle\prod_{i=1}^r \frac{\sigma_i\left(\beta_i\right)}{\beta_i}.
	\]
\end{Pf}

We remark that $(\cdot,\cdot)_{L/K}$ induces 
\[
\begin{array}{cccc}
	U_1					&\rightarrow	&\mathrm{Hom}(\widehat{G}, \mathbb{Q}/\mathbb{Z})	& (\cong G)	\\
	\rotatebox{90}{$\in$}		&			&\rotatebox{90}{$\in$}							&		\\[-4pt]
	u					&\mapsto		&(u, \cdot)_{L/K}								&		
\end{array}
\]
under the assumption in Proposition \ref{bilinear}. The isomorphism is Pontryagin dual.

\begin{Rmk}
	The above map (named $P_{L/K}$) is indicated by the followings. We regard $U_1$ as a group of $k$-rational points
	of pro-algebraic group $\mathbb{U}_1$ over $k$, where $k$ is the residue field of $K$. Let $\pi_0(\mathbb{U}_1)$ be
	the pro-finite group of pro-connected components in $\mathbb{U}_1$. 
	By Serre-Hazewinkel's geometric local class field theory (see \cite{Ha}), we have the natural isomorphism
\[
  \pi_0(\mathbb{U}_1) \cong G.
\]
	Then $P_{L/K}$ corresponds with composition of the above isomorphism and the natural map $U_1 \rightarrow \pi_0(\mathbb{U}_1)$.
\end{Rmk}

\subsection{Coleman power series and $(\cdot ,\cdot )_{L/K}$}\label{CTCPS}
Now we start to prove Theorem \ref{maintheorem1}.  Let $K$ be as above, and take a uniformizer $\pi$. Let $\mathcal{F}$ be a Lubin-Tate module of $\pi$ over $\mathcal{O}_K$, and fix a generator $\xi = (\xi_n)_n \in \varprojlim_n \mathcal{F}(n)$ of Tate module.  Let $H$ be a finite unramified extension of $K$. Put
$H_n := H(\xi_n)$, $K_n := K(\xi_n)$, $H_\infty := \underset{n}\cup H_n$ , $G := \Gal{H_\infty}{H}$. Let $\mathrm{rec}_H$ be the reciprocity map given by the class field theory
\[
	\mathrm{rec}_H : \mathcal{O}_H^\times \xrightarrow{\sim} G.
\]
Let $\mathrm{Col}_u$ denotes the Coleman power series for $u = (u_n)_n \in \varprojlim_n{\mathcal{O}_{H_n}^\times}$
and $L$ be a finite extension of $H$ contained in $H_\infty$. Put $u_L := N_{H_n/L}(u_n)$ where $n$ is large enough.
Then the following theorem holds.
 
\begin{Thm}\label{cnj}
	Assume that $u_{H} = 1$. Then for any character $\chi : G \rightarrow \mathbb{Q}/\mathbb{Z}$ whose order is finite, the following equality holds. 
	\[
		\mathrm{rec}_{H}^{-1}((u_{H_\chi}, \cdot)_{H_\chi/H}) = \mathrm{Col}_u(0)
	\]
	where $H_\chi = (H_\infty)^{\mathrm{Ker}\chi}$. 
\end{Thm}

\begin{Pf}
	Let $b_u \in \mathcal{O}_H^\times$ denote an element which satisfies
	$\mathrm{rec}_{H}^{-1}((u_{H_\chi}, \cdot)_{H_\chi/H}) = b_u$. We have to prove that $b_u = \mathrm{Col}_u(0)$.
	Put $u' = (u_n')_n := (\mathrm{rec}_{H}(b_u)(\xi_n)/\xi_n)_n \in \varprojlim_n \mathcal{O}_{H_n}^\times$ and
	$u'' = (u''_n)_n := (u_n/u'_n)_n$. By definition, we have $\mathrm{rec}_{H}^{-1}((u'_{H_\chi}, \cdot)_{H_\chi/H}) = b_u$ and
	$\mathrm{rec}_{H}^{-1}((u''_{H_\chi}, \cdot)_{H_\chi/H}) = 1$. Therefore $(u''_{H_\chi}, \cdot)_{H_\chi/H} = 0$.
	We remark that the sequence $(\xi_n)_n$ is norm coherent since the definition of $[\pi]$, 
	so the sequence $(u_n')_n$ is. 
	Then we have
	\[
		\mathrm{Col}_{u'}(T) = \displaystyle\frac{\mathrm{Col}_{\mathrm{rec}_{H}(b_u)(\xi)}(T)}{\mathrm{Col}_{(\xi)}(T)}
		= \displaystyle\frac{b_uT + (\text{higher degree terms})}{T}.
	\]
	It implies $\mathrm{Col}_{u'}(0) = b_u$. Therefore we have 
	\[
	\mathrm{Col}_{u''}(0) = \mathrm{Col}_{u}(0)/\mathrm{Col}_{u'}(0) = \mathrm{Col}_{u}(0)/b_u.
	\]
	Thus we have to prove $\mathrm{Col}_{u''}(0) = 1$.
	Fix an integer $n \geq 0$. 
	Let $\mathcal{F}[\pi^{n}]$ denote the finite flat group scheme of $\pi^{n}$-torsion points of $\mathcal{F}$
	over $\mathcal{O}_K$ and $R_{n}$ the coordinate ring of $\mathcal{F}[\pi^{n}]$. 
	Then $R_{n}$ is local and finite flat $\mathcal{O}_K$-algebra whose rank is 
	equal to the $n$-th power of the order of the residue field of $K$.
	We regard $\xi_n$ as a rational point of $\mathcal{O}_{K_n}$ at $\mathcal{F}(n)$ and then we have 
	\[
		\mathrm{Spec}\,\mathcal{O}_{K_n} \rightarrow \mathcal{F}(n).
	\]
	This is an open immersion and the image is equal to $\mathcal{F}(n)\backslash\mathcal{F}(n-1)$. So there exists a surjection
	\[
		r_n : R_n \twoheadrightarrow \mathcal{O}_{K_n}.
	\]
	Since $R_n$ is local, $r_n$ induces $R_n^\times \twoheadrightarrow \mathcal{O}_{K_n}^\times$. We also denote it by $r_n$.
	Let $r'_n$ denote the surjection $(R_n\otimes_{\mathcal{O}_K} \mathcal{O}_{H})^\times \twoheadrightarrow
	\mathcal{O}_{H_n}^\times$ induced by $r_n$.
	Let $\widetilde{u_n} \in (R_n\otimes_{\mathcal{O}_K}\mathcal{O}_{H})^\times$ denote a lifting of $u''_n$ which satisfies
	$r'_n(\widetilde{u_n}) = u''_n$.
	It is sufficient that the following equality holds independently of a choice of $\widetilde{u_n}$:
	\begin{align}\label{eq3}
		N_{H_n/H}(\widetilde{u_n})
		\equiv \mathrm{Col}_{u''}(0) \quad \mathrm{mod}\,\pi^{n+1}\mathcal{O}_{H}
	\end{align}
	If it holds for any $n$, we have $N_{H_n/H}(\widetilde{u_n}) = \mathrm{Col}_{u''}(0)$.
	Applying Lemma \ref{14} to $u''$, there exist a finite unramified extension $H'$, an integer $r$,
	elements $\beta_1,\ldots,\beta_r \in \mathcal{O}_{H'H_n}^\times$ and $\sigma_1,\ldots,\sigma_r \in \Gal{H'H_n}{H'}$ such that 
	$u''_n = \displaystyle\prod_{i=1}^r \frac{\sigma_i\left(\beta_i\right)}{\beta_i}$.
	Let $\widetilde{\beta_i} \in (R_n\otimes_{\mathcal{O}_K}\mathcal{O}_{H'})^\times$ denotes a lifting of $\beta_i$ 
	which satisfies $r'_n(\widetilde{\beta_i}) = \beta_i$. Then we have
	\[
		\widetilde{u_n} = \displaystyle\prod_{i=1}^r \frac{\sigma_i\left(\widetilde{\beta_i}\right)}{\widetilde{\beta_i}}.
	\]
	Since it is clear that $N_{H'H_n/H'}\circ\sigma_i = N_{H'H_n/H'}$ for all
	$\sigma_i \in \Gal{H'H_n}{H'}$. Therefore we have
	\[
		N_{H'H_n/H'}\left( \displaystyle\frac{\sigma_i\left(\widetilde{\beta_i}\right)}{\widetilde{\beta_i}}\right) = 1
	\]
	for all $0 \leq i \leq r$. It implies $N_{H_n/H}(\widetilde{u_n}) = 1$.
	At last we prove the claim (\ref{eq3}). To prove the independence of a choice of $\widetilde{u_n}$, we now show that
	\[
	\text{if}\ x \in \mathrm{Ker}\, r_n\ \text{then}\ N_{H_n/H}(x) \equiv 1 \quad \mathrm{mod}\, \pi^{n+1}\mathcal{O}_K.
	\]
	Let $x \in R_n\otimes_{\mathcal{O}_K}\mathcal{O}_{H}$ be an element which satisfies $r_n(x) = 1$ for any $n \geq 0$,
	$t_n$ a surjection from $R_n\otimes_{\mathcal{O}_K}\mathcal{O}_{H}$ to $R_{n-1}\otimes_{\mathcal{O}_K}\mathcal{O}_{H}$
	induced by $\mathcal{F}[\pi^{n}] \hookrightarrow \mathcal{F}[\pi^{n+1}]$ where $R_{-1} := \mathcal{O}_K$.
	Since $R_n \cong \mathcal{O}_K[[T]]/\left( [\pi^{n+1}](T) \right)$, we have
	\[
	(R_n\otimes_{\mathcal{O}_K}\mathcal{O}_{H})^\times\cong\left(\mathcal{O}_{H}[[T]]/\left( [\pi^{n+1}](T) \right)\right)^\times.
	\]
	Since $x$ satisfies $r_n(x) = 1$, there exists a polynomial $g(T) \in \mathcal{O}_{H}[T]$ such that
	\[
	x-1 \equiv h(T)g(T) \quad \mathrm{mod} [\pi^{n+1}](T)
	\]
	where $h(T) := [\pi^{n+1}](T)/[\pi^n](T)$. It is clear that there exists a polynomial $h_1(T) \in \mathcal{O}_{H}[T]$ such that $h(T) \equiv \pi h_1(T)
	\ \mathrm{mod}\, [\pi^n](T)$ and then $t_n(h(T)g(T)) \equiv 0 \quad \mathrm{mod}\ \pi R_{n-1}\otimes_{\mathcal{O}_K}\mathcal{O}_{H}$.
	Hence we have:
	\[
	t_n(x) \equiv 1 \quad \mathrm{mod}\ \pi R_{n-1}\otimes_{\mathcal{O}_K}\mathcal{O}_{H}.
	\]
	Next it is also clear that $r_{n-1}(N_{H_n/H_{n-1}}(x)) = N_{{K_n'}/K_{n-1}'}(r_n(x)) = 1$.
	Moreover if $x$ satisfies $r_n(x) = 1$ and $t_n(x) \equiv 1 \quad \mathrm{mod}\, \pi R_{n-1}\otimes_{\mathcal{O}_K}\mathcal{O}_{H}$ then we have
	\[
	t_{n-1}(N_{H_n/H_{n-1}}(x)) \equiv 1 \quad \mathrm{mod}\, \pi^2 R_{n-2}\otimes_{\mathcal{O}_K}\mathcal{O}_{H}.
	\]
	By repeating the above discussion, we have
	\[
	t_0(N_{H_n/H_{0}}(x)) \equiv 1 \quad \mathrm{mod}\, \pi^{n+1} R_{0}\otimes_{\mathcal{O}_K}\mathcal{O}_{H}
	\] 
	We remark that $N_{H_n/H_{n-1}}(y) = (t_0(y))(N_{H_0/H}(r_n(y)))$
	for any $y \in R_0\otimes_{\mathcal{O}_K}\mathcal{O}_{H}$. Hence we have 
	\[
	N_{H_n/H}(x) \equiv 1 \quad \mathrm{mod}\, \pi^{n+1}\mathcal{O}_K.
	\]
	Let $f_n$ be the image of $\mathrm{Col}_{u''}(T)$ under the natural projection $\mathcal{O}_{H}[[T]]^\times
	\twoheadrightarrow \left(\mathcal{O}_{H}[[T]]/\left( [\pi^{n+1}](T) \right)\right)^\times$. Therefore we have
	\[
		N_{H_n/H}(\widetilde{u_n})
		\equiv N_{H_n/H}(f_n) \quad \mathrm{mod}\,\pi^{n+1}\mathcal{O}_{H}
	\]
	and then we have to prove $N_{H_n/H}(f_n) = \mathrm{Col}_{u''}(0)$ in $H$.
	For $0 \leq i \leq n$, let $\xi_i$ denote the following immersion
	\[
		\mathrm{Spec}\,\mathcal{O}_{K_i} \underset{\text{open}}\hookrightarrow \mathcal{F}(i)
		\underset{\text{closed}}\hookrightarrow \mathcal{F}(n).
	\] 
	Since $\mathcal{F}(n+1)$ has the unit point $\mathrm{Spec}\,\mathcal{O}_K \hookrightarrow \mathcal{F}(n)$,
	we have following maps:
	\begin{align*}
		R_n &\rightarrow \mathcal{O}_{K_i} (\text{ for all } 0 \leq i \leq n) \\
		R_n &\rightarrow \mathcal{O}_K.
	\end{align*}
	Here it is well known that the map
	\[
		R_n\otimes_{\mathcal{O}_K} H \rightarrow H \times \displaystyle\prod_{i=0}^n H_i
	\]
	induced by the map $R_n \rightarrow \mathcal{O}_{K}\times \displaystyle\prod_{i=0}^n \mathcal{O}_{K_i}$ is isomorphism.
	Then the image of $f_n$ under the map $R_n\otimes_{\mathcal{O}_K} H \rightarrow H_i$ is equal to $\mathrm{Col}_{u''}(\xi_i)$ and
	under the map $R_n\otimes_{\mathcal{O}_K} H \rightarrow H$ is equal to $\mathrm{Col}_{u''}(0)$. Therefore we have
	\begin{align*}
		N_{H_n/H}(f_n)
			&= N_{H/H}\mathrm{Col}_{u''}(0) \times \displaystyle\prod_{i=0}^n N_{H_i/H}\mathrm{Col}_{u''}(\xi_i) \\
			&= \mathrm{Col}_{u''}(0) \times \displaystyle\prod_{i=0}^n N_{H_i/H}(u''_i) \\
			&= \mathrm{Col}_{u''}(0) \times \displaystyle\prod_{i=0}^n N_{H_i/H}(u_i) \\
			&= \mathrm{Col}_{u''}(0) \times \displaystyle\prod_{i=0}^n u_{H}\\
			&= \mathrm{Col}_{u''}(0).
	\end{align*}
	Since $u_{H} = 1$ and this completes the proof.
\end{Pf}

\begin{Cor}\label{20}
	Under the assumption in Theorem \ref{cnj}, we have
	\[
		(u_{H_\chi} , \chi)_{H_\chi/H} = \chi \circ \mathrm{rec}_{H}(\mathrm{Col}_u(0))
	\]
	for any $u \in \varprojlim_n{\mathcal{O}_{H_n}^\times}$ and $\chi \in \widehat{G}$ 
\end{Cor}
\begin{Pf}
	It follows immediately by Theorem \ref{cnj}.  
\end{Pf}

	Let $u \in \varprojlim_n \mathcal{O}_{H_n}^\times$ be an element which satisfies $u_{K} = 1$ and $\chi : G \rightarrow \mathbb{Q}/\mathbb{Z}$ be a character of finite order $d$.
	Assume that there exist a positive integer $d$ and $a \in \mathcal{O}_{H}^\times$ such that $u_{H} = a^d$ and $H_{\chi}$ 
	has no nontrivial $d$-th roots of unity.
	Then, since $N_{H/K}(a) =1$ by using Hilbert 90, there exists an element $b \in \mathcal{O}_{H}^\times$
	such that $a = \mathrm{Frob}(b)/b$ where $\mathrm{Frob} \in \Gal{H}{K}$ is an arithmetic Frobenius.
	Then Theorem \ref{cnj}
	shows the following. This corollary is used to prove Theorem \ref{fsc} mainly.

\begin{Cor}\label{16}
	\[
		(u_{H_\chi}/a , \chi)_{H_\chi/H} = \chi \circ \mathrm{rec}(\mathrm{Col}_u(0) / b^d)
	\]
\end{Cor}

\begin{Pf}
	Now $N_{H/K}(a) = 1$ so there exists an element $\widetilde{u} \in \varprojlim_n \mathcal{O}_{H_n}^\times$ such that
	$\widetilde{u}_{H} = a$. Via the Coleman power series $\mathrm{Col}_{\widetilde{u}}$ for $\widetilde{u}$,
	it is easy to show that $\mathrm{Frob}(\mathrm{Col}_{\widetilde{u}}(0))/\mathrm{Col}_{\widetilde{u}}(0) = \widetilde{u}_{H} = a$, 
	and there exists $c \in \mathcal{O}_K^\times$ such that $\mathrm{Col}_{\widetilde{u}}(0) = bc$. Since 
	$N_{H_\chi/K}(a/\widetilde{u}_{H_\chi}^d) = (a / \widetilde{u}_{H})^d = 1$, we have
	\[
		(a/\widetilde{u}_{H_\chi}^d , \chi)_{H_\chi/H} = 0
	\]
	by Proposition \ref{symbol0}.	Using Corollary \ref{20}, we have
	\begin{align*}
		(u_{H_\chi}/a , \chi)_{H_\chi/H} &= \left((u_{H_\chi}/\widetilde{u}_{H_\chi}^d)(\widetilde{u}_{H_\chi}^d/a) , \chi \right)_{H_\chi/H}
				 \\
				&= (u_{H_\chi}/\widetilde{u}_{H_\chi}^d , \chi)_{H_\chi/H} - (a/\widetilde{u}_{H_\chi}^d , \chi)_{H_\chi/H} \\
				&= \chi \circ \mathrm{rec}(\mathrm{Col}_{u}(0) / \mathrm{Col}_{\widetilde{u}}(0)^d) \\
				&= \chi \circ \mathrm{rec}(\mathrm{Col}_{u}(0) / b^d c^d).
	\end{align*}
	Since $d$ is the order of $\chi$, the last term is equal to $\chi \circ \mathrm{rec}(\mathrm{Col}_u(0) / b^d)$.
\end{Pf}


\section{Drinfeld module and Siegel units}\label{start}
Now we start to construct our Euler systems. In this section we prepare Drinfeld modules and Siegel units.
This contents in this section is based on \cite{KY}.

\subsection{Drinfeld module}
Let $p$ be an odd prime number and $q$ a power of $p$, $\mathbb{F}_q$ the finite field of $q$ elements 
$C$ a nonsingular projective geometrically irreducible curve over $\mathbb{F}_q$ and $F$ the function field of $C$.
We regard the closed points of $C$ as the primes of $F$ and let $F_v$, where $v$ is closed point of $C$, be the completion of 
$F$ at $v$, ${\cal O}_v$ the integer ring of $F_v$, $\pi_v \in {\cal O}_v$ a prime element, $k(v)$ the residue field of ${\cal O}_v$.
Fix a closed point $\infty$ of $C$. Let $p_\infty$ be the order of $k(\infty)$. Since
$C \backslash \left\{ \infty \right\}$ is affine, we can identify it with 
$\mathrm{Spec}\, A$ where $A = \Gamma(C\backslash \left\{ \infty \right\} ,  {\cal O}_C )$. 
From now on, unless otherwise stated, all primes of $F$ are assumed to be different from $\infty$.
Next we define Drinfeld modules and Drinfeld modular varieties \cite{Dr}.

\begin{Def}[Drinfeld module]
	Let $S$ be an $A$-scheme. {\bf A Drinfeld module} of rank 1 over $S$ is an scheme $E$ in $A$-modules
	over $S$ which satisfies following conditions:
	\begin{enumerate}
		\item Zariski locally on $S$, the scheme $E$ is isomorphic to the additive group scheme $\mathbb{G}_{a, S}$ as a commutative group scheme.
		\item We denote the $A$-action on $E$ by $\mathrm{act} : A \rightarrow \mathrm{End}_{S-gp}(E)$.
			For every $a \in A \setminus \{ 0 \}$, the morphism $\mathrm{act}(a) : E \rightarrow E$ is 
			finite, locally free of constant degree $|a|_{\infty}$ where $|a|_{\infty}$ is absolute value at $\infty$.
		\item The $A$-action on Lie $E$ induced by $\mathrm{act}$ coincides with the $A$-action on Lie $E$ 
			which comes from the structure homomorphism $E \rightarrow S$.  
	\end{enumerate}
\end{Def}

When we say a Drinfeld module, we assume that its rank is equal to 1. 

\begin{Def}[Drinfeld Modular Variety]
	Let $N$ be a nonzero torsion $A$-module.
	Put $U_{N} := \mathrm{Spec}\ A \setminus \mathrm{Supp}\ N$. Let $S$ be a $U_{N}$-scheme and $E =  (E, \mathrm{act})$
	a {\rm Drinfeld module} over $S$. A {\bf level $N$ structure} on $E$ is a monomorphism
	\begin{equation*}
		\mathrm{lev} : N_S \hookrightarrow E
	\end{equation*}
	of schemes in $A$-modules over $S$, where $N_S$ is a constant group scheme. We consider the following 
	contravariant functor from the category of $U_N$-schemes to the category of sets:
		\[\begin{array}{ccc}
		(U_N\text{-schemes}) 	&\rightarrow	& (\text{Sets})	\\
		\rotatebox{90}{$\in$}		&			&\rotatebox{90}{$\in$}	\\[-4pt]
		S						&\mapsto	&(\text{the set of isomorphism classes of }(E, \mathrm{act}, \mathrm{lev}))
	\end{array}\]
	If $N \neq 0$, there exists an $U_N$-scheme functor $\mathcal{M}_N$ such that the functor is isomorphic to the functor which makes $S$ 
	correspond with $\mathrm{Hom}_{U_N}(S ,  \mathcal{M}_N)$. The $U_N$-scheme isomorphic is unique
	up to canonical isomorphisms due to
	Yoneda's lemma. We call $\mathcal{M}_N$ the {\bf Drinfeld Modular Variety} of rank 1 with level $N$-structures.
 \end{Def}

For two torsion $A$-modules $N, N'$ and an embedding $N \hookrightarrow N'$, let $r_{N, N'}$ denote the following morphism from $\mathcal{M}_{N'}$ 
to $\mathcal{M}_{N} \underset{U_{N}}{\times} U_{N'}$:
\[
\begin{array}{rccc}
	r_{N, N'} :	&\mathcal{M}_{N'} 				&\rightarrow	&\mathcal{M}_{N} \underset{U_{N}}{\times} U_{N'}	\\
			&\rotatebox{90}{$\in$}			&			&\rotatebox{90}{$\in$}							\\[-4pt]
			&(E, \mathrm{act}, \mathrm{lev})	&\mapsto	&(E, \mathrm{act}, \mathrm{lev} |_N)
\end{array}
\]
where $\mathrm{lev} |_N$ is restriction of $\mathrm{lev}$ on $N$.
Similarly, for two torsion $A$-modules $N, N''$ and a surjection $N \twoheadrightarrow N''$, we call it $\pi$ for only this section whose kernel is of finite length, let $m_{N, N''}$ denote the following morphism from $\mathcal{M}_{N}$ to 
$\mathcal{M}_{N''} \underset{U_{N''}}{\times} U_{N}$:
\[
\begin{array}{rccc}
	m_{N, N''} :	&\mathcal{M}_{N}				&\rightarrow	&\mathcal{M}_{N''} \underset{U_N''}{\times} U_N		\\
			&\rotatebox{90}{$\in$}			&			&\rotatebox{90}{$\in$}							\\[-4pt]
			&(E, \mathrm{act}, \mathrm{lev})	&\mapsto	&(E'', \mathrm{act}'', \mathrm{lev}'')
\end{array}
\]
where $E''=E/\mathrm{lev}(\mathrm{Ker}\ \pi)$ and $\mathrm{act}'', \mathrm{lev}''$ are the maps induced by $\mathrm{act}$ and $\mathrm{lev}$ respectively.

\subsection{Theta function}
We define the theta function in this section. Let $S$ be a reduced scheme over $\mathrm{Spec}\, A$ and $(E, \mathrm{act})$ a
Drinfeld module on $S$. We regard $S$ as a closed subscheme in $E$ via the zero section $S \hookrightarrow E$.

\begin{Prop} \label{existf}
	For any nonzero element $a \in A$, let $N_a$ denote the norm map 
	\[
	\Gamma(E \backslash \mathrm{Ker}\, \mathrm{act}(a) ,  {\mathcal O}_{E}^\times )
	\rightarrow \Gamma(E \backslash S ,  {\mathcal O}_{E}^\times )
	\]
	with respect to a finite flat morphism 
	\[
	\mathrm{act} (a):E\backslash \mathrm{Ker}\, \mathrm{act}(a) \rightarrow E \backslash S.
	\]
	Then there exists the function $f \in \Gamma(E \backslash S ,  {\cal O}_{E}^\times )$
	unique up to $\mu_{p_\infty -1} (S)$ which satisfies the following:
	\begin{enumerate}
		\item $N_a(f) = f$ for any nonzero element $a \in A$.
		\item The order $\mathrm{ord}_S (f)$ of $f$ in $S$ is equal to $p_\infty -1$. 
	   \end{enumerate}
\end{Prop}
\begin{Pf}
	See \cite[Lemma2.1]{KY}.
 \end{Pf}

From the above, the function $f^{p_\infty -1} \in \Gamma(E \backslash S ,  {\cal O}_{E}^\times )$ is uniquely defined.

\begin{Def}[Theta function]
	For a Drinfeld module $(E, \mathrm{act})$ , let $f$ be an element satisfying the condition in Proposition \ref{existf}. Let 
	\[
	\vartheta_{E/S} := f^{p_\infty -1}.
	\]
	We call $\vartheta_{E/S}$ the {\bf theta function} of $(E, \mathrm{act})$.
 \end{Def}

\begin{Rmk}\label{theta1}
	Let $S, S'$ be schemes over Spec\, $A$, $\mathrm{mor}:S' \rightarrow S$ a morphism from $S'$ to $S$, 
	$\mathrm{mor}_E:E \underset{S}{\times} S' \rightarrow E$ base change with respect to $\mathrm{mor}$.
	Then we have
	\[\mathrm{mor}_E^\ast\, \vartheta_{E/S} = \vartheta_{E \underset{S}{\times} S'/S'}.  \]
\end{Rmk}

\begin{Rmk}\label{theta2}
	Let $E, E'$ be Drinfeld modules on $S$ and $\mathrm{iso}:E \rightarrow E'$ be an isogeny.
	Let $N_\mathrm{iso}$ denote the norm map corresponding to $\mathrm{iso}$. Then
	\[N_\mathrm{iso}\, \vartheta_{E/S} = \vartheta_{E'/S}. \]
\end{Rmk}

\subsection{Siegel unit}\label{defsie}
Let $N$ be a nonzero torsion $A$-module as above, $E_N \rightarrow \mathcal{M}_N$ the universal Drinfeld module,
and $\mathrm{lev}:N_{\mathcal{M}_N} \hookrightarrow E_N$ the universal level structure.
Let us apply the argument in the last paragraph to the case when $S =  \mathcal{M}_{N}$. Via the zero section we regard $\mathcal{M}_{N}$ as a closed subscheme of $E_N$. 
We remark that $\mathcal{M}_{N}$ is smooth over $U_N$, in particular $\mathcal{M}_{N}$ is reduced.
Therefore Proposition \ref{existf} assures the existence of a theta function 
$\vartheta_{E_N / \mathcal{M}_N} \in \Gamma(E_N \backslash \mathcal{M}_N ,  {\cal O}_{E_N \backslash \mathcal{M}_N}^\times )$.

\begin{Def}[Siegel unit]
	For an element $b$ in $N \backslash \{ 0 \}$, let $\mathrm{lev}_b : \mathcal{M}_N \rightarrow E_N$ denote the restriction of
	$\mathrm{lev}$ to $\mathcal{M}_{N}$, where we regard $\mathcal{M}_N = \{ b \} \times \mathcal{M}_N$ as a subscheme of  
	$N_{\mathcal{M}_{N}} = \displaystyle \coprod_{b \in N} \mathcal{M}_N$. Put
	\[
	\sie{N, b} := \mathrm{lev}_b^\ast\, \vartheta_{E_N / \mathcal{M}_N} \in 
	\Gamma(\mathcal{M}_N, {\cal O}_{\mathcal{M}_N}^\times )
	\]
	and we call it a {\bf Siegel unit}.
 \end{Def}

Let $N'$ be an $A$-module of finite length, $N$ an $A$-submodule of $N'$. Then both $\mathcal{M}_N ,  \mathcal{M}_{N'}$ are not empty. Let $N''$ be an $A$-module of finite length and 
$\alpha : N \twoheadrightarrow N''$ a surjection. Then we can prove the following proposition using 
Fact \ref{theta1}, \ref{theta2}.

\begin{Prop}[distribution relation]\label{dist}
	For any $b \in N \backslash \{ 0 \} ,  b'' \in N'' \backslash \{ 0 \}$, following equalities hold.
	\begin{align*}
		r_{N', N}^\ast\, \sie{N, b} &= \sie{N', b} \\
		m_{N, N''}^\ast\, \sie{N'', b''} &= \prod_{b \in N ,  \alpha(b)=b''} \sie{N, b} 
	\end{align*}
\end{Prop}

\begin{Prop}\label{unit}
	If there exist two different maximal ideals of $A$ which divide $\mathrm{Ann}_A(b) := \{ a\in A | ab = 0\}$, then
	$\sie{N, b}$ is a unit of the integral closure of $A$ in $\Gamma(\mathcal{M}_N ,  {\cal O}_{\mathcal{M}_N})$.
 \end{Prop}

Let $I \neq A$ be an ideal of $A$ and $v$ a prime of $F$ not dividing $I$. Consider the case $N = A/v^nI$. For any integers
$n \geq 0$, we can construct Siegel units $\sie{A/v^nI ,  \{1\}}$ as above.
For simplicity, we write $m_n$ for $m_{A/v^nI ,  A/v^{n-1}I}$.

\begin{Prop}\label{hecke}
	The following equality between Siegel units holds
	\begin{align*}
		m_{n,  \ast}\, \sie{A/v^nI ,  \{1\}} = \left\{
		\begin{array}{ll}
			\sie{A/v^{n-1}I ,  \{1\} } &(\text{ if } n\geq 2), \\
			(1-T_v)\sie{A/I ,  \{1\} } &(\text{ if } n = 1) \\
		\end{array}
		\right. 
	\end{align*}
	where $T_v$ is a Hecke operator determined by $v$ and satisfies $r_{A/I ,  A/vI}^\ast = T_v^\ast m_1^\ast$. 
\end{Prop}
\begin{Pf}
	Using the fact that $m_n^\ast$ is injective, it is enough to prove that 
	\[
	m_n^\ast m_{n,  \ast}\, \sie{A/v^nI ,  \{1\}} = \left\{
		\begin{array}{ll}
			m_n^\ast \sie{A/v^{n-1}I ,  \{1\} } &(\text{ if } n\geq 2), \\
			m_n^\ast (1-T_v)\sie{A/I ,  \{1\} } &(\text{ if } n = 1). \\
		\end{array}
		\right.
	 \]
	For $n \geq 1$, let 
	\begin{align*}
		\pi_n &: A/v^nI \twoheadrightarrow A/v^{n-1}I , \\
		\pi_n^\times &: (A/v^nI)^\times \twoheadrightarrow (A/v^{n-1}I)^\times
	\end{align*}
	be natural surjections. There exists a natural isogeny
	\[
		\mathrm{iso} : E_{A/v^nI} \rightarrow m_n^\ast\, E_{A/v^{n-1}I}.
	\]
	Recall that $N_{\mathrm{iso}}$ denote the norm map with respect to $\mathrm{iso}$.
	By the definition of $\vartheta$ and Proposition \ref{dist}, we obtain
	\begin{align}
	\label{eq1} N_\mathrm{iso} \vartheta_{E_{A/v^nI}/ \mathcal{M}_{A/v^nI}} 
		= \vartheta_{m_n^\ast\ E_{A/v^{n}I} / \mathcal{M}_{A/v^{n}I}}.
	\end{align}
	The following diagram is commutative.
	\begin{align*}
	\begin{xy}
		(0, 0)*{E_{A/v^nI}}="e", 
		(40, 0)*{m_n^\ast\ E_{A/v^nI}}="m", 
		(20, 7)*{\circlearrowright}, 
		(20, 15)*{\mathcal{M}_{A/v^{n}I}}="mm", 
		\ar_{\mathrm{iso}} "e";"m", 
		\ar^{m_n^\ast\, \mathrm{lev}_1} "mm";"m", 
		\ar_{\mathrm{lev}_1} "mm";"e", 
	\end{xy}
	\end{align*}
	Apply $\mathrm{lev}_1^\ast\,  \mathrm{iso}^\ast = (m_n \mathrm{lev}_1)^\ast $ to the both sides of formula (\ref{eq1}) and 
	calculate, we have
	\begin{align*}
		\text{right hand side} &= (m_n \mathrm{lev}_1)^\ast\, 
			\vartheta_{m_n^\ast\ E_{A/v^{n}I}/ \mathcal{M}_{A/v^{n}I}} \\
		&= m_n^\ast\ \mathrm{lev}_1^\ast\, 
			\vartheta_{E_{A/v^{n-1}I}/ \mathcal{M}_{A/v^{n-1}I}} \\
		&= m_n^\ast\ \sie{A/v^{n-1}I ,  \{1\}}.
	\end{align*}
	On the other hand, we have
	\begin{align*}
		\text{left hand side} &= \mathrm{lev}_1^\ast\ \mathrm{iso}^\ast\ N_{\mathrm{iso}}\ 
			\vartheta_{E_{A/v^{n}I}/ \mathcal{M}_{A/v^{n}I}} \\
		&= \mathrm{lev}_1^\ast\ \mathrm{iso}^\ast\ \mathrm{iso}_\ast\ 
			\vartheta_{E_{A/v^{n}I}/ \mathcal{M}_{A/v^{n}I}} \\
		&= \mathrm{lev}_1^\ast (\displaystyle\prod_{a \in \mathrm{Ker}\, \pi_n}t_a^\ast\ 
			\vartheta_{E_{A/v^{n}I}/ \mathcal{M}_{A/v^{n}I}}) \\
		&= \displaystyle\prod_{a \in \mathrm{Ker}\, \pi_n}(\mathrm{lev}_1^\ast\ t_a^\ast\ 
			\vartheta_{E_{A/v^{n}I}/ \mathcal{M}_{A/v^{n}I}}) \\
		&= \left\{
		\begin{array}{ll}
			\displaystyle\prod_{a \in (\mathrm{Ker}\, \pi_n)^\times}(\mathrm{lev}_a^\ast\ 
				\vartheta_{E_{A/v^{n}I}/ \mathcal{M}_{A/v^{n}I}}) & \\ & (\text{ if } n\geq 2 ) \\
			\displaystyle\prod_{a \in (\mathrm{Ker}\, \pi_n)^\times}(\mathrm{lev}_a^\ast\ 
				\vartheta_{E_{A/v^{n}I}/ \mathcal{M}_{A/v^{n}I}}) &
				\times (\mathrm{lev}_b^\ast\, \vartheta_{E_{A/v^{n}I}/ \mathcal{M}_{A/v^{n}I}}) \\ & (\text{ if } n=1 ) \\
		\end{array}\right. \\
		&= \left\{
		\begin{array}{ll}
			m_n^\ast\ m_{n ,  \ast}\ \sie{A/v^nI ,  \{1\}} &(\text{ if } n\geq 2 ) \\
			(\mathrm{lev}_b^\ast\, \vartheta_{E_{A/v^{n}I}/ \mathcal{M}_{A/v^{n}I}})\times
			m_n^\ast\ m_{n ,  \ast}\ \sie{A/v^nI ,  \{1\}} &(\text{ if } n=1 ) \\
		\end{array}\right. 	
	\end{align*}
	where $t_a$ denotes the action of the element $a$ in $\mathrm{Ker}\, \pi$, i.e.
	\begin{align*}
		t_a : E_{A/v^nI} =& E_{A/v^nI} \underset{\mathcal{M}_{A/v^nI}}{\times} \mathcal{M}_{A/v^nI} \\
		\xrightarrow{\mathrm{id}\times\mathrm{lev}_a} & E_{A/v^nI} \underset{\mathcal{M}_{A/v^nI}}{\times} E_{A/v^nI}
		\xrightarrow{\text{group law}} E_{A/v^nI}
	\end{align*}
	and $b$ is the element of $A/vI$ which is mapped to $(1,0)$ under the natural isomorphism $A/vI \cong A/I \times A/v$. Remark that we have used the fact that
	\begin{align*}
		\mathrm{Ker}\, \pi = \left\{
		\begin{array}{llll}
			(v^{n-1}A/v^nA) &\cong (1+v^{n-1}A/v^nA)  &= \mathrm{Ker}\, \pi^\times &(\text{ if } n\geq 2 ) \\
			A/v &\cong (A/v)^\times \cup \{ b \} &= \mathrm{Ker}\, \pi^\times \cup \{ b \} &(\text{ if } n=1 ). \\
		\end{array}
		\right. 
	\end{align*}
	Then if $n\geq2$, we have
	\[
	m_n^\ast\ \sie{A/v^{n-1}I ,  \{1\}} =
	m_n^\ast\ m_{n ,  \ast}\,  \sie{A/v^nI ,  \{1\}}
	\]
	as desired. If $n=1$, there is an extra term $\mathrm{lev}_b^\ast\, 
	\vartheta_{E_{A/v^{n}I}/ \mathcal{M}_{A/v^{n}I}}$. Now we see that it can be written in terms of the Hecke operator $T_v$.
	Here $I$ and $v$ are relatively prime so we have $A/vI \cong A/I\oplus A/v$. Via the natural embedding 
	$A/I \hookrightarrow A/vI$ we regard $A/I$ as a direct summand of $A/vI$. Let us consider the morphism
	\begin{align*}
	\begin{array}{rccc}
	r_{A/I, A/vI} :	&\mathcal{M}_{A/vI} 		&\rightarrow	&\mathcal{M}_{A/I} \underset{U_{A/I}}{\times} U_{A/vI}	\\
			&\rotatebox{90}{$\in$}			&			&\rotatebox{90}{$\in$}							\\[-4pt]
			&(E, \mathrm{act}, \mathrm{lev})	&\mapsto	&(E, \mathrm{act}, \mathrm{lev} |_{A/I}).
	\end{array}
	\end{align*}
	From its definition, we obtain the isomorphism
	\[
	E_{A/vI} \xrightarrow{\sim} E_{A/vI} \underset{\mathcal{M}_{A/I}}{\times} \mathcal{M}_{A/vI}.
	\]
	Moreover, since $b \in A/vI$ is equal to the image of $1 \in A/I$ under the embedding $A/I \hookrightarrow A/vI$,
	we can prove that the composite of $\mathrm{lev}_b$ and the above isomorphism is equal to, 
	\[
	\mathrm{lev}_1\times \mathrm{id} : \mathcal{M}_{A/vI} = \mathcal{M}_{A/I} \underset{\mathcal{M}_{A/I}}{\times} \mathcal{M}_{A/vI}
	\hookrightarrow E_{A/I} \underset{\mathcal{M}_{A/I}}{\times} \mathcal{M}_{A/vI}.
	\]
	Since the Hecke operator $T_v$ satisfies the equality $r_{A/I ,  A/vI}^\ast =  m_1^\ast T_v^\ast$, we have
	\begin{align*}
	\mathrm{lev}_b^\ast\,  \vartheta_{E_{A/vI}/ \mathcal{M}_{A/vI}} &= m_1^\ast T_v^\ast  
		\mathrm{lev}_1 \vartheta_{E_{A/I}/ \mathcal{M}_{A/I}} \\
	&= m_1^\ast T_v^\ast \sie{A/I ,  \{1\}}.
	\end{align*}
	Then put in order,
	\begin{align*}
		m_1^\ast\ \sie{A/I ,  \{1\}} &=
			m_1^\ast T_v^\ast \sie{A/I ,  \{1\}} \times m_1^\ast\ m_{1 ,  \ast}\ \sie{A/vI ,  \{1\}} \\
		&= m_1^\ast(T_v^\ast \sie{A/I ,  \{1\}} \times m_{1 ,  \ast}\ \sie{A/vI ,  \{1\}})
	\end{align*}
	therefore we have
	\[ 
	m_1^\ast  (1-T_v^\ast ) \sie{A/vI ,  \{1\}}
		= m_1^\ast\ m_{1 ,  \ast}\ \sie{A/vI ,  \{1\}}. 
	\]
\end{Pf}

In the next section, we write $\mathcal{M}_{N, F}$ fors $\mathcal{M}_N \underset{U_N}{\times} {\rm Spec}\,  F$ for simplicity.


\section{Euler Systems}
In this section we prepare a tower of field extensions analogues to the cyclotomic extension $\mathbb{Q}(\mu_{p^n})$
in algebraic number fields. More preciously, we consider $\mathcal{O} ( \mathcal{M}_{A/v^nI, F} )$ and
$\Gamma(\mathcal{M}_{A/v^nI, F})$  for any $n \geq 0$.
And then we confirm that Siegel units belong to the tower and satisfy the norm relations of Euler systems. 

\subsection{Notation}\label{seceu}
Let $p, F, \infty, A$ be as in the first paragraph of Section 2. Let $I$ be an ideal of $A$ not equal to $A$.
Put $\hat{A} :=\displaystyle \varprojlim_{0 \neq J: \text{ideal of}\, A} A/J$. Let $v$ be a prime of $F$ not dividing $I$,
$\mathbb{A}_F$ the ring of adeles of $F$, and $\mathbb{A}_F^\times$ the group of ideles of $F$.

We briefly mention what is $\mathcal{O} ( \mathcal{M}_{A/v^nI, F} )$. First it is proved by Drinfeld that 
$\mathcal{M}_{A/v^nI, F}$ equals to the spectrum of finite abelian extension of $F$ \cite[\S 8-Theorem1]{Dr}.
Hence $\mathcal{M}_{A/v^nI, F}$ is, as a topological space, consists of a single element and
the structure sheaf is a finite abelian extension of $F$. Moreover $\infty$ splits completely in $\mathcal{O}(\mathcal{M}_{A/v^n I,F})$. 
We regard the tower of extensions $\mathcal{O} ( \mathcal{M}_{A/v^nI, F} )$ as an analogue of $\mathbb{Q}(\mu_{p^n})$ in algebraic field.
Next, it is also proved by Drinfeld that the reciprocity map of the global class field theory induces an isomorphism
\begin{align}\label{equal1}
	\Gal{\fd{v^n}}{F} \cong F^\times F_\infty^\times \backslash \mathbb{A}_F^\times /
		\left(\hat{A}^\times \cap (1+v^nI\hat{A})\right).
\end{align}

Let $F_{v}$ be the completion of $F$ at $v$ and $\mathcal{O}_{v}$ the integer ring. Similarly let $\mathcal{O} ( \mathcal{M}_{A/I, F} )_w$
denote the completion of $\mathcal{O} ( \mathcal{M}_{A/I, F} )$ at $w$ and $\mathcal{O}_w$ the integer ring.
Now let us consider the extension $\mathcal{O} ( \mathcal{M}_{A/I, F} )$ over $F$. Let $v$ be a prime of $F$ not dividing $I$ which split
completely in $\mathcal{O} ( \mathcal{M}_{A/I, F} ) / F$. Let $M \in \mathbb{N}$ be a power of $p_\infty$. Let $\psi_{v}$ be
a continuous surjective homomorphism
\[
\psi_{v} : \mathcal{O}_{v}^\times \twoheadrightarrow \mathbb{Z}/M\mathbb{Z}
\]
which factors through $(\mathcal{O}_{v}/\mathfrak{p}_{v}^{n_v})$ for some $n_v$, 
where $\mathfrak{p}_{v}$ is the maximal ideal of the local ring $\mathcal{O}_{v}$.

\begin{Def}
	Let $\Psi_M$ denote the set of finite sets of pairs of a prime $v$ and a homomorphism $\psi_{v}$  satisfying the conditions as above, i.e. 
	\begin{eqnarray*}
		\Psi_M := \left\{ \{(v_1, \psi_{v_1}) ,  ... ,  (v_r, \psi_{v_r}) \} \left|
		\begin{array}{ll}
			r \geq 0 ,  v_i \ (i = 1, ..., r) \text{ is a prime of } F  \\ 
			\text{which is prime to } I, \\
			\text{split completely in } \mathcal{O} ( \mathcal{M}_{A/I, F} ) / F,\\
			\text{and both } v_i \text{ and } v_j \ (i \neq j)	\text{ are prime, }\\
			\psi_{v_i} \text{ is a continuous surjective} \\ 
			\text{homomorphism as above. }	
		\end{array}
		\right. \right\}
	\end{eqnarray*}
	where $r=0$ means $\emptyset \in \Psi_M$.
 \end{Def}
For any $v_i$, fix an integer $n_i := n_{v_i}$ depending on $v_i$ as above. Put
\[
(\bfvec{v}_r ,  \psi_{\bfvec{v}_r}) := \{(v_1, \psi_{v_1}) ,  ... ,  (v_r, \psi_{v_r}) \}
\]
for simplicity. Let ``$(v, \psi_v) \in \Psi_M$'' denote ``$\{ (v, \psi_v) \} \in \Psi_M$''.

\subsection{Norm equation}
For an element $(\bfvec{v}_r ,  \psi_{\bfvec{v}_r}) \in \Psi_M$, let us consider the field extension corresponding to it:
\begin{align*}
	\begin{xy}
		(0, 0)*{F.}="ff0", 
		(0, 10)*{\fd{}}="ff1", 
		(0, 20)*{\fd{v_1^{n_1}}} ="ff2", 
		(0, 30)*{\vdots} ="ff3", 
		(0, 38)*{\fd{v_1^{n_1}\cdots v_r^{n_r}}} ="ff4", 
		\ar @{-} (0, 3);(0, 7)
		\ar @{-} (0, 13);(0, 17)
		\ar @{-} (0, 23);(0, 27)
		\ar @{-} (0, 31);(0, 35)		
	\end{xy}
\end{align*}

We have the following.
\begin{Prop}\label{gal}
	$\Gal{\fd{v_1^{n_1}\cdots v_r^{n_r}}}{\fd{}} \cong 
	(A /v_1^{n_1}\cdots v_r^{n_r}A)^\times \cong \displaystyle \prod_{i=1}^r (A /v_i^{n_i}A)^\times $
 \end{Prop}

\begin{Pf}
	The latter isomorphism is clear by Chinese reminder theorem since any $v_i$ is prime to $v_j (i \neq j)$.
	Let us prove the former. Only in this proof, we put 
	$v := v_1^{n_1}\cdots v_r^{n_r}$ for simplicity. The above diagram and the Galois group (\ref{equal1}) show that 
	the following sequence of abelian groups is exact:
	\begin{align*}
			0 \rightarrow \Gal{\fd{v}}{\fd{}}
			&\rightarrow F^\times F_\infty^\times \backslash \mathbb{A}_F^\times /(\hat{A}^\times \cap (1+vI\hat{A})) \\
			&\rightarrow F^\times F_\infty^\times \backslash \mathbb{A}_F^\times /(\hat{A}^\times \cap (1+I\hat{A}))
			 \rightarrow 0.
	\end{align*}
	So we have to prove
	\begin{align*}
		\mathrm{Ker}\left(F^\times F_\infty^\times \backslash \mathbb{A}_F^\times /(\hat{A}^\times \cap (1+vI\hat{A}))
		\rightarrow  F^\times F_\infty^\times \backslash \mathbb{A}_F^\times /(\hat{A}^\times \cap (1+I\hat{A}))\right)
			\\ \cong (A/vA)^\times.
	\end{align*}
	First, by the definition of it, we have 
	\begin{align*}
		\mathrm{Ker}&\left(F^\times F_\infty^\times \backslash \mathbb{A}_F^\times /	(\hat{A}^\times \cap (1+I\hat{A}))
			\rightarrow F^\times F_\infty^\times \backslash \mathbb{A}_F^\times /(\hat{A}^\times \cap (1+I\hat{A}))\right) \\
		&= \hat{A}^\times \cap (1+I\hat{A})/\left( (\hat{A}^\times\cap(1+I\hat{A})) \cap F^\times F_\infty^\times (
			\hat{A}^\times \cap (1+vI\hat{A}))\right) \\
		&= \hat{A}^\times \cap (1+I\hat{A})/\left( F^\times F_\infty^\times \cap (\hat{A}^\times\cap(1+I\hat{A})\cdot
			\hat{A}^\times \cap (1+vI\hat{A}))\right).
	\end{align*}
	Here it is clear that
	\begin{align*}
		F^\times F_\infty^\times \cap (\hat{A}^\times\cap(1+I\hat{A}))
		&\cong F^\times F_\infty^\times \cap \hat{A}^\times \\
		&\cong F^\times \cap (F_\infty^\times\cdot\hat{A}^\times) \\
		&= \{ x \in F^\times \mid \text{ for any }v\neq\infty, x\in\mathcal{O}_{v} \} \\
		&= \Gamma(C\backslash{\infty}, \mathcal{O}_C^\times) = A^\times,
	\end{align*}
	so recall that $I\neq A$, we have
	\begin{align*}
		\hat{A}^\times \cap (1+I\hat{A})&/\left( F^\times F_\infty^\times \cap (\hat{A}^\times\cap(1+I\hat{A})\cdot
			\hat{A}^\times \cap (1+vI\hat{A}))\right) \\
		&\cong \hat{A}^\times\cap(\hat{A}^\times \cap (1+I\hat{A}))\cdot(\hat{A}^\times \cap (1+vI\hat{A})) \\
		&\cong 1\cdot(\hat{A}^\times \cap (1+vI\hat{A})).
	\end{align*}
	Therefore,
	\begin{align*}
		\mathrm{Ker}&\left(F^\times F_\infty^\times \backslash \mathbb{A}_F^\times/(
			\hat{A}^\times \cap (1+vI\hat{A}))\right)	\rightarrow
			(F^\times F_\infty^\times \backslash \mathbb{A}_F^\times)/(\hat{A}^\times\cap(1+I\hat{A}))\\
		&\cong (\hat{A}^\times \cap (1+I\hat{A}))/\hat{A}^\times \cap (1+vI\hat{A}) \\
		&\cong (1+I\hat{A})/(1+vI\hat{A}) \\
		&\cong (A/vA)^\times. 
	\end{align*}
\end{Pf}

Take $\psi_{v_i}$ apart as follows: 
\[
\mathcal{O}_{F_{v_i}}^\times \rightarrow (\mathcal{O}_{F_{v_i}}/\mathfrak{p}_{v_i}^{n_i})^\times 
\twoheadrightarrow \mathbb{Z}/M\mathbb{Z}
\]
and put
\[
\overline{\psi_{v_i}} : (\mathcal{O}_{F_{v_i}}/\mathfrak{p}_{v_i}^{n_i})^\times \rightarrow \mathbb{Z}/M\mathbb{Z}.
\]
the latter map. 
Let $E(\bfvec{v}_r ,  \psi_{\bfvec{v}_r})$ be a subfield of $\fd{v_1^{n_1}\cdots v_r^{n_r}}$ over $\fd{}$ whose the Galois group 
$\Gal{\fd{v_1^{n_1}\cdots v_r^{n_r}}}{E(\bfvec{v}_r ,  \psi_{\bfvec{v}_r})}$ is equal to
${\rm Ker}\, (\overline{\psi_{v_1}} \times \cdots \times \overline{\psi_{v_r}})$.
For the sake of convenience if $r=0$, we define $E(\bfvec{v}_0 ,  \psi_{\bfvec{v}_0}) := \fd{}$. It is clear that
\[
	\Gal{E(\bfvec{v}_r ,  \psi_{\bfvec{v}_r})}{\fd{}} \cong \displaystyle \prod_{i=1}^r (\mathbb{Z}/M\mathbb{Z}).
\]
Put
\[
G_i := \Gal{E(\bfvec{v}_i ,  \psi_{\bfvec{v}_i})}{E(\bfvec{v}_{i-1} ,  \psi_{\bfvec{v}_{i-1}})}
\]
then there exists a natural isomorphism 
$\displaystyle \prod_{i=1}^r G_i \cong \Gal{E(\bfvec{v}_r ,  \psi_{\bfvec{v}_r})}{\fd{}}$.
Let $\sie{A/v_1^{n_1} \cdots v_r^{n_r} I ,  \{ 1 \}} \in \fd{v_1^{n_1} \cdots v_r^{n_r}}$ be a Siegel unit as section \ref{defsie}. Put 
\[
	\sie{(\bfvec{v}_r ,  \psi_{\bfvec{v}_r})} := {\rm Norm}_{\fd{v_1^{n_1} \cdots v_r^{n_r}}/E(\bfvec{v}_r ,  \psi_{\bfvec{v}_r})}\, 
	\sie{A/v_1^{n_1} \cdots v_r^{n_r} I ,  \{ 1 \}}.
\]
For any $G_i$ let us choose an element $\sigma_i \in G_i$ whose image under the isomorphism $G_i \cong \mathbb{Z}/M\mathbb{Z}$
is equal to 1. We define operators $\bfvec{N}_i$ and $\bfvec{D}_i$ as 
\begin{eqnarray*}
	\bfvec{N}_i &:=& \sum_{\tau \in G_i} \tau \in \mathbb{Z}[G_i] \\
	\bfvec{D}_i &:=& \sum_{j=1}^{\sharp G_i-1} j\sigma_i^j \in \mathbb{Z}[G_i]
\end{eqnarray*}
where we write an action for multiplicative group $\fd{v_1^{n_1} \cdots v_r^{n_r}}^\times$ additively. We will often use this additive notation. 

\begin{Prop}\label{nd}
	The following equality holds for any $i$.
	\[ (\sigma_i - 1)\bfvec{D}_i = M - \bfvec{N}_i  \]
\end{Prop}

\begin{Pf}
	Calculating both sides with taking care of $\sharp G_i = M$, the proposition is clear.
 \end{Pf}

Put $\bfvec{N} := \displaystyle \prod_{i=1}^r \bfvec{N}_i ,  \bfvec{D} := \displaystyle \prod_{i=1}^r \bfvec{D}_i$.
Then the Siegel unit $\sie{(\bfvec{v}_r ,  \psi_{\bfvec{v}_r})}$ satisfies the following conditions.

\begin{Thm}[Norm relations]\label{norm}
For any $(\bfvec{v}_r ,  \psi_{\bfvec{v}_r})$ and $(v_i ,  \psi_{v_i}) \in (\bfvec{v}_r ,  \psi_{\bfvec{v}_r})$, 
	\begin{description}
		\item[ES1] $\sie{(\bfvec{v}_r ,  \psi_{\bfvec{v}_r})} \in  (E(\bfvec{v}_r ,  \psi_{\bfvec{v}_r}))^\times$.
		\item[ES2] $\sie{(\bfvec{v}_r ,  \psi_{\bfvec{v}_r})}$ is a unit of the integral closure.
		\item[ES3] $\bfvec{N}_i\, \sie{(\bfvec{v}_r ,  \psi_{\bfvec{v}_r})} = 
					(1-{\rm Frob}_i)\sie{(\bfvec{v}_{r \backslash i} ,  \psi_{\bfvec{v}_{r \backslash i}})}$.
	\end{description}
Where ${\rm Frob}_i$ is the Frobenius of $v_i$ in $G_r$, $(\bfvec{v}_{r \backslash i} ,  \psi_{\bfvec{v}_{r \backslash i}})$
is an element that $(v_i ,  \psi_{v_i})$ is removed from $(\bfvec{v}_r ,  \psi_{\bfvec{v}_r})$.
 \end{Thm}

\begin{Pf}
	Drinfeld \cite[$\S$8, Thm.1, Cor.]{Dr} proved that $\mathcal{M}_{A/v^nI, F}$ is a spectrum of a
	finite abelian extension over the function field $F$. We have the following dictionary: 
	\begin{align*}
		\mathcal{M}_{A/I, F} &\longleftrightarrow \text{ the spectrum of a finite abelian extension over } F \\
		\mathcal{M}_{A/I} &\longleftrightarrow \text{ normalization of } U_{A/I} \text{ in } \mathcal{M}_{A/I, F}   \\
		m_{A/v_i^nI ,  A/v_i^{n-1}I ,  \ast} &\longleftrightarrow \text{ the norm operator } \\
		T_{v_i}^\ast &\longleftrightarrow \text{ the } \mathrm{Frob}_{i} 
	\end{align*}
	for a prime $v_i$.
	This leads {\bf ES2} using Proposition \ref{unit}. {\bf ES1} is clear from the definition.
	Remark that Proposition \ref{hecke} can be written as 
	\[
	\bfvec{N}_i \sie{A/v_i^nI ,  \{1\}} = \left\{
	\begin{array}{ll}
		\sie{A/v_i^{n-1}I ,  \{1\} } &(\text{ if } n\geq 2), \\
		(1-\mathrm{Frob}_i)\sie{A/I ,  \{1\} } &(\text{ if } n = 1). \\
	\end{array}
	\right.
	\]
	So take $I$ as $(v_1^{n_1}\cdots v_r^{n_r}/v_i^{n_i})I$, we have
	\[
	\bfvec{N}_i \sie{A/v_i^{n_i}(v_1^{n_1}\cdots v_r^{n_r}/v_i^{n_i})I ,  \{1\}} = (1-{\rm Frob}_i)\sie{A/I ,  \{1\}}
	\]
	Therefore {\bf ES3} is hold.
\end{Pf}


\section{Finite-singular comparison equalities}\label{kol}
In this section, we construct Kolyvagin's derivative class starting from our Euler system and formulate the
``finite-singular comparison equalities'' (Theorem \ref{fsc} and \ref{cheb}). Our assumption is that $M$ is power of $p$.
Note that if $M$ is prime to $p$, similar relations are proved by H.Oukhaba and S.Vigui$\acute{{\rm e}}$ in \cite{OS}.

\subsection{Derivative class}\label{kap}
For Theorem \ref{norm}, the following lemma holds.

\begin{Lem}\label{3lem}
	Let $(\bfvec{v}_r ,  \psi_{\bfvec{v}_r}) \in \Psi_M$. Then
	\[
	\bfvec{D} \sie{(\bfvec{v}_r ,  \psi_{\bfvec{v}_r})}\in
	(E(\bfvec{v}_r ,  \psi_{\bfvec{v}_r})^\times/(E(\bfvec{v}_r, \psi_{\bfvec{v}_r})^\times)^M)^{G}
	\]
	where $G = \displaystyle \prod_{i=1}^r G_i$.
 \end{Lem}

\begin{Pf}
	We prove this by induction on the cardinality of the finite set $r$ included in $(\bfvec{v}_r ,  \psi_{\bfvec{v}_r})$.
	Let $\sigma_i$ be a fixed generator of $G_i \quad (1\leq i \leq r)$. When $r=1$, i.e.
	$(\bfvec{v}_r ,  \psi_{\bfvec{v}_r}) = (v_1, \psi_{v_1})$, by the Proposition \ref{nd} and {\bf ES3} in \ref{norm},
	\begin{eqnarray*}
		(1 - \sigma_1)\bfvec{D}_1\sie{(v_1 ,  \psi_{v_1})} &=& (M - \bfvec{N}_1)\sie{(v_1 ,  \psi_{v_1})} \\
		& \equiv & -\bfvec{N}_1\sie{(v_1 ,  \psi_{v_1})}\quad {\rm mod}\ (E(\bfvec{v}_1, \psi_{\bfvec{v}_1})^\times)^M \\
		& = & -(1-{\rm Frob}_{1})\sie{A/I ,  \{1\}}\quad {\rm mod}\ (E(\bfvec{v}_1, \psi_{\bfvec{v}_1})^\times)^M.
	\end{eqnarray*}
	Here ${\rm Frob}_{1}$ acts trivially on $\fd{}$ and $\sie{A/I ,  \{1\}} \in \fd{}$, so the last term is $0$.
	
	Let us assume for all $r-1$, the claim holds. In other words, for all $(\bfvec{v}_{r-1} ,  \psi_{\bfvec{v}_{r-1}}) \in \Psi_M$ and
	for all $\sigma \in \displaystyle\prod_{i=1}^{r-1} G_i$, the following holds:
	\begin{eqnarray*}
		(1 - \sigma)\displaystyle\prod_{i=1}^{r-1}\bfvec{D}_i\sie{(\bfvec{v}_{r-1} ,  \psi_{\bfvec{v}_{r-1}})} 
			\equiv 0 \quad {\rm mod}\ (E(\bfvec{v}_{r-1}, \psi_{\bfvec{v}_{r-1}})^\times)^M
	\end{eqnarray*}
	Now we prove that
	\[
		(1- \sigma)\displaystyle\prod_{i=1}^{r}\bfvec{D}_i \sie{(\bfvec{v}_{r} ,  \psi_{\bfvec{v}_{r}})} 
		\equiv 0\quad {\rm mod}\ (E(\bfvec{v}_1, \psi_{\bfvec{v}_1})^\times)^M 
	\]
	for $(\bfvec{v}_r , \psi_{\bfvec{v}_r}) \in \Psi_M$ and $\sigma \in G = \displaystyle\prod_{i=1}^{r}G_i$.
	Since $G = \displaystyle\prod_{i=1}^{r}G_i$ is generated by $\sigma_i \in G_i (i = 1,\ldots,r)$, we can assume that
	$\sigma = \sigma_i$. Furthermore, by the symmetry, we may assume $\sigma = \sigma_r$ without loss of generality. Then we have
	\begin{eqnarray*}
		(1 - \sigma_r)\displaystyle\prod_{i=1}^{r}\bfvec{D}_i\sie{(\bfvec{v}_{r}, \psi_{\bfvec{v}_{r}})} 
				&=& (M - \bfvec{N}_r)\displaystyle\prod_{i=1}^{r-1}\bfvec{D}_i\sie{(\bfvec{v}_{r}, \psi_{\bfvec{v}_{r}})} \\
		&\equiv & -\bfvec{N}_r\displaystyle\prod_{i=1}^{r-1}\bfvec{D}_i\sie{(\bfvec{v}_{r}, \psi_{\bfvec{v}_{r}})}
				\quad {\rm mod}\ (E(\bfvec{v}_{r}, \psi_{\bfvec{v}_{r}})^\times)^M \\
		&=& -(1-{\rm Frob}_{r})\displaystyle\prod_{i=1}^{r-1}\bfvec{D}_i\sie{(\bfvec{v}_{r-1}, \psi_{\bfvec{v}_{r-1}})}
				\quad \\ && {\rm mod}\ (E(\bfvec{v}_{r}, \psi_{\bfvec{v}_{r}})^\times)^M.
	\end{eqnarray*}
	Due to our assumption, the last term is equal to $0$.
 \end{Pf}

\begin{Lem}
	For any $(\bfvec{v}_r ,  \psi_{\bfvec{v}_r}) \in \Psi_M$, there exists an element
	\[
	\kappa_{(\bfvec{v}_r ,  \psi_{\bfvec{v}_r})} \in \fd{}^\times/(\fd{}^\times)^M
	\]
	such that
	\[
	\kappa_{(\bfvec{v}_r ,  \psi_{\bfvec{v}_r})} \equiv \bfvec{D}\sie{(\bfvec{v}_r ,  \psi_{\bfvec{v}_r})}
		\ {\rm mod}\ (E(\bfvec{v}_r ,  \psi_{\bfvec{v}_r})^\times)^M.  
	\]
\end{Lem}

\begin{Pf}
	Recall that $M$ is power of $p$, there is no $M$-th roots of unity in $E(\bfvec{v}_{r} ,  \psi_{\bfvec{v}_{r}})$. 
	By taking the Galois cohomology of the following short exact sequence:
	\begin{align*}
		&0 \rightarrow E(\bfvec{v}_{r} ,  \psi_{\bfvec{v}_{r}})^\times \xrightarrow{\text{ power of }M} E(\bfvec{v}_{r} ,  \psi_{\bfvec{v}_{r}})^\times \\
		&\quad \xrightarrow{\text{ projection }} E(\bfvec{v}_{r} ,  \psi_{\bfvec{v}_{r}})^\times / (E(\bfvec{v}_{r} ,  \psi_{\bfvec{v}_{r}})^\times)^M \rightarrow 0.
	\end{align*}
	We have the following exact sequence.
	\begin{align*}
		&0 \rightarrow \fd{}^\times \xrightarrow{\text{ power of } M} \fd{}^\times \\
		&\quad \rightarrow \left(E(\bfvec{v}_{r} ,  \psi_{\bfvec{v}_{r}})^\times / (E(\bfvec{v}_{r} ,  \psi_{\bfvec{v}_{r}})^\times)^M\right)^G
		\rightarrow 0.
	\end{align*}
	The exactness of the last term follows from Hilbert 90. Hence the inclusion of $\fd{}^\times$ into $E(\bfvec{v}_{r} ,  \psi_{\bfvec{v}_{r}})^\times$
	leads a natural isomorphism
	\begin{align*}
		\fd{}^\times / (\fd{}^\times)^M \cong
		\left(E(\bfvec{v}_{r} ,  \psi_{\bfvec{v}_{r}})^\times / (E(\bfvec{v}_{r} ,  \psi_{\bfvec{v}_{r}})^\times)^M\right)^G.
	\end{align*}
	Lemma \ref{3lem} says $\bfvec{D} \sie{(\bfvec{v}_r ,  \psi_{\bfvec{v}_r})}
	\in(E(\bfvec{v}_r ,  \psi_{\bfvec{v}_r})^\times/(E(\bfvec{v}_r, \psi_{\bfvec{v}_r})^\times)^M)^{G}$,
	put $\kappa_{(\bfvec{v}_r ,  \psi_{\bfvec{v}_r})}$ be an element corresponding to 
	$\bfvec{D} \sie{(\bfvec{v}_r ,  \psi_{\bfvec{v}_r})}$ by the isomorphism above.
\end{Pf}

We call the family $\{ \kappa_{(\bfvec{v}_r ,  \psi_{\bfvec{v}_r})} \}_r$ of the above element
$\kappa_{(\bfvec{v}_r ,  \psi_{\bfvec{v}_r})}$ {\bf Kolyvagin's derivative classes}.

\subsection{Finite-singular comparison equalities}\label{FSCE}
Let $\mathscr{I} := \underset{\lambda}\bigoplus \mathbb{Z} \lambda$ be the group of fractional ideals of $\fd{}$, 
where $\lambda$ ranges over the primes of $\fd{}$. For any prime $v$ of $F$, let us identify $v$ as a prime ideal of $A$, 
and define $\mathscr{I}_v := \underset{\lambda \mid v}\bigoplus \mathbb{Z} \lambda$.
It is clear that $\mathscr{I} = \underset{v}\bigoplus \mathscr{I}_v$.
Let $(y) \in \mathscr{I}$ denote the  principal ideal generated by $y$ where $y$ is an element of $\fd{}^\times$, 
$(y)_v \in \mathscr{I}_v$ the $v$-part of $(y)$, $[y] \in \mathscr{I}/M\mathscr{I}$ the projection of $(y)$,
and $[y]_v \in \mathscr{I}_v / M\mathscr{I}_v$ the $v$-part of $[y]$.
Then the following condition holds.

\begin{Thm}[finite-singular comparison equalities]\label{fsc}
	Let $(\bfvec{v}_r ,  \psi_{\bfvec{v}_r}) \in \Psi_M$, $(v, \psi_v) \in \Psi_M$. Then
	\begin{align*}
		\text{ For any } i = 1, ..., r , &v \neq v_i\quad \\ &\Rightarrow\quad [\kappa_{(\bfvec{v}_r ,  \psi_{\bfvec{v}_r})}]_v = 0 \\
		\text{ For some } i = 1, ..., r \text{ such that } &v = v_i\quad \\ &
		\Rightarrow\quad [\kappa_{(\bfvec{v}_r ,  \psi_{\bfvec{v}_r})}]_v
		 = \underset{\lambda \mid v}\sum \psi_v 
			(\iota_\lambda(\kappa_{(\bfvec{v}_{r \backslash i},\psi_{\bfvec{v}_{r \backslash i}})}))\lambda 
	\end{align*}
	where $\iota_\lambda$ is the isomorphism $\mathcal{O}_\lambda \xrightarrow{\sim} \mathcal{O}_{F_v}$. 
\end{Thm}

This theorem is one of the main result in this paper. These equalities are necessary to prove the Iwasawa main conjecture. In classical case, $F=\mathbb{Q}$,
these equalities (and proof of it) are more simple since we can construct a map $\psi_v$ which is not depend on $v$. However this method could not apply in our case because the characteristic of $F$ is not zero. We prove these equalities using Theorem \ref{cnj}.

\subsection{Proof of Theorem \ref{fsc}}
Let $v_i \left( 1 \leq i \leq r \right)$ be a prime of $F$ which does not divide the ideal $I , \infty$ and split completely in
$\fd{}/F$. Fix a prime $w_i$ of $\fd{}$ above $v_i$. Let $w_i^j \left( 1 \leq j \leq r \right)$ be the prime
of $\fd{v_1^{n_1} \cdots v_j^{n_j}}$ above $v_i$ and $\mathcal{O}_{v_i}$ be the integral ring 
of the completion of $F$ at $v_i$, $\mathcal{O}_{w_i^j}$ the integer ring of the completion of $\fd{}$ at $w_i^j$. 
Remark that $\mathcal{O}_{w_i^j} \cong \mathcal{O}_{v_i}$ since $v_i$ splits completely in $\fd{}/F$. 
We consider a Drinfeld module $E_{A/I} \cong \mathrm{Spec}\ R[T]$
whose the zero section $\mathcal{M}_{A/I} \rightarrow E_{A/I}$ satisfies $T \mapsto 0$.
Make $E_{A/I}$ base change from $\mathcal{M}_{A/I}$ to $\mathcal{O}_{w_i^j}$, and put
\[
	E_{\mathcal{O}_{w_i^j}} := E_{A/I} \underset{\mathcal{M}_{A/I}}\times \mathrm{Spec}\ \mathcal{O}_{w_i^j}.
\]
Remark that $E_{\mathcal{O}_{w_i^j}} \cong \mathrm{Spec}\ \mathcal{O}_{w_i^j}[T]$. The following map is induced by group law of 
$E_{\mathcal{O}_{w_i^j}} \times E_{\mathcal{O}_{w_i^j}} \rightarrow E_{\mathcal{O}_{w_i^j}}$
	\begin{align*}
	\begin{xy}
		(0, 0)*{\mathcal{O}_{w_i^j}[T]}, 
		(30, -1)*{\mathcal{O}_{w_i^j}[T] \underset{\mathcal{O}_{w_i^j}}\otimes \mathcal{O}_{w_i^j}[T]}, 
		(50, 0)*{\cong}, 
		(65, 0)*{\mathcal{O}_{w_i^j}[T,S]}, 
		(0, -5)*{\rotatebox{90}{$\in$}}, 
		(65, -5)*{\rotatebox{90}{$\in$}}, 
		(0, -8)*{T}, 
		(65, -8)*{P(T,S).}, 
		\ar @{->} (7,0);(11,0), 
		\ar @{|->} (4,-8);(55,-8), 
	\end{xy}
	\end{align*}
We regard $P(T,S)$ as a formal power series by $\mathcal{O}_{w_i^j}[T,S] \subseteq \mathcal{O}_{w_i^j}[[T,S]]$ and
$\mathcal{F}$ a formal group law over $\mathcal{O}_{w_i^j}$,
since it satisfies the definition of formal group law by the zero section satisfying $T \mapsto 0$.  
By the definition of Drinfeld module, for a uniformizer $\pi$ of $\mathcal{O}_{v_i}$, it is clear that
\[
	[\pi]_F(T) \equiv T^{q_{w_i^j}} \ \mathrm{mod}\ \pi
\]
where $q_{w_i^j}$ is the order of the residue field of $\mathcal{O}_{w_i^j}$. 
So $\mathcal{F}$ is the Lubin-Tate module of $\pi$. Let $\mathcal{F}(n)$ be the set of $\pi^{n+1}$-th roots of unity, 
$\xi := (\xi_n)_n$ a generator of Tate module $\varprojlim_n \mathcal{F}(n)$. Put
$\mathcal{O}_{w_i^j}(n) := \mathcal{O}_{w_i^j}[\mathcal{F}(n)]$. Then we have
\[
	\mathcal{M}_{A/v^nI} \cong \mathrm{Spec}\ \mathcal{O}_{w_i^j}(n).
\]
Recall from section \ref{start},
	\begin{align*}
	\begin{xy}
		(0, 0)*{(E , \mathrm{act} , \mathrm{lev})}, 
		(0, 4)*{\rotatebox{90}{$\in$}}, 
		(0, 8)*{\mathcal{M}_{A/v_1^{n_1}\cdots v_j^{n_j} I}}, 
		(0, 20)*{E_{A/v_1^{n_1}\cdots v_j^{n_j}I}}, 
		(70, 0)*{(E/\mathrm{lev}(\mathrm{Ker}(A/v_1^{n_1}\cdots v_j^{n_j}I \twoheadrightarrow A/v_1^{n_1}\cdots v_j^{n_j -1}I))) , \mathrm{act}' , \mathrm{lev}')}, 
		(70, 4)*{\rotatebox{90}{$\in$}}, 
		(70, 8)*{\mathcal{M}_{A/v_1^{n_1}\cdots v_j^{n_j -1}I}}, 
		(70, 20)*{E_{A/v_1^{n_1}\cdots v_j^{n_j -1}I}}, 
		\ar @{->} (0,17);(0,11), 
		\ar @{->} (70,17);(70,11), 
		\ar @{->}^{m_{A/v_1^{n_1}\cdots v_j^{n_j}I , A/v_1^{n_1}\cdots v_j^{n_j -1}I}} (12,8);(55,8), 
		\ar @{|->} (12,0);(20,0), 
	\end{xy}
	\end{align*}
Let $m_{n_j}$ denote the map repeated above for $n_j$ times
\[
	m_{n_j} : \mathcal{M}_{A/v_1^{n_1}\cdots v_j^{n_j}I} \rightarrow \mathcal{M}_{A/v_1^{n_1}\cdots v_{j-1}^{n_{j-1}}I}.
\]
Here $v_j$ splits completely in $\fd{}/F$, in particular $v_j$ is a principle ideal. So there exists an element $\pi_{v_j} \in A$ such that
\[
	v = \pi_{v_j}A.
\]
Let $(\pi_v) : E/v_1^{n_1}\cdots v_j^{n_j}I \rightarrow E/v_1^{n_1}\cdots v_j^{n_j}I$ be the map given by the multiplication by $\pi_{v_j}$. Then we have
\[
	m_{n_j}^* E_{A/I} \cong E_{A/v_1^{n_1}\cdots v_j^{n_j}I} / \mathrm{lev}(\mathrm{Ker}({\pi_{v_j}}^{n_j})) \cong E_{A/v_1^{n_1}\cdots v_j^{n_j}I}
\]
where the last isomorphism is induced by $\pi^{n_j}$. Let $\phi_{n_j}$ be the composition
\begin{align*}
	&(A/v_1^{n_1}\cdots v_{j-1}^{n_j -1}I \times \pi_{v_j}^{-n_j}A/A)_{\mathcal{M}_{A/v_1^{n_1}\cdots v_j^{n_j}I}} \\ \cong 
	&(A/v_1^{n_1}\cdots v_{j-1}^{n_j -1}I \times A/v_1^{n_1}\cdots v_j^{n_j}A)_{\mathcal{M}_{A/v_1^{n_1}\cdots v_j^{n_j}I}} \\ 
	\xrightarrow{\mathrm{lev}} &E_{A/v_1^{n_1}\cdots v_j^{n_j}I} \cong m_{n_j}^* E_{A/v_1^{n_1}\cdots v_{j-1}^{n_j -1}I}
\end{align*}
and $\phi_{n_j}'$ the composition
\begin{align*}
	&(A/v_1^{n_1}\cdots v_{j-1}^{n_j -1}I)_{\mathcal{M}_{A/v_1^{n_1}\cdots v_j^{n_j}I}} 
	\xrightarrow{\mathrm{lev}} &E_{A/v_1^{n_1}\cdots v_j^{n_j}I} \cong m_{n_j}^* E_{A/v_1^{n_1}\cdots v_{j-1}^{n_j -1}I}.
\end{align*}
Let us embed $\mathcal{M}_{A/v_1^{n_1}\cdots v_j^{n_j}I}$ into
\[
(A/v_1^{n_1}\cdots v_{j-1}^{n_j -1}I \times \pi_v^{-n}A/A)_{\mathcal{M}_{A/v_1^{n_1}\cdots v_{j}^{n_j}I}} = \underset{b \in (A/v_1^{n_1}\cdots v_{j-1}^{n_j -1}I \times \pi_v^{-n}A/A)}\coprod \mathcal{M}_{A/v_1^{n_1}\cdots v_{j}^{n_j}I}.
\]
We now calculate what a Siegel unit will be. 
\begin{align*}
	\sie{A/v_1^{n_1}\cdots v_{j}^{n_j}I , \{ 1 \}} 
	&= \mathrm{lev}^* \vartheta_{E_{A/v_1^{n_1}\cdots v_{j}^{n_j}I}/\mathcal{M}_{A/v_1^{n_1}\cdots v_{j}^{n_j}I}} \\
	&= \phi_{n_j}(\pi_{v_j}^{-n_j} , \pi_{v_j}^{-n_j})^* m_{n_j}^* \vartheta_{E_{A/v_1^{n_1}\cdots v_{j}^{n_j}I}/\mathcal{M}_{A/v_1^{n_1}\cdots v_{j}^{n_j}I}} \\
	&= \phi_{n_j}((\pi_{v_j}^{-n_j} , 0) + (0 ,\pi_{v_j}^{-n_j}))^* m_{n_j}^* \vartheta_{E_{A/v_1^{n_1}\cdots v_{j}^{n_j}I}/\mathcal{M}_{A/v_1^{n_1}\cdots v_{j}^{n_j}I}} \\
	&= \phi_{n_j}(0 , \pi_{v_j}^{-n_j})^* m_{n_j}^* t_{\phi_{n_j}'(\pi_{v_j}^{-n_j})} \vartheta_{E_{A/v_1^{n_1}\cdots v_{j}^{n_j}I}/\mathcal{M}_{A/v_1^{n_1}\cdots v_{j}^{n_j}I}}
\end{align*}
where $t_{\phi_{n_j}'(\pi_{v_j}^{-n_j})}$ appearing the last term is the translation by $\phi_{n_j}'(\pi_{v_j}^{-n_j})$,
\begin{align*}
	t_{\phi_{n_j}'(\pi_{v_j}^{-n_j})} \ :& \ E_{A/v_1^{n_1}\cdots v_{j}^{n_j}I} \xrightarrow{\sim} E_{A/v_1^{n_1}\cdots v_{j}^{n_j}I} \underset{\mathcal{M}_{A/v_1^{n_1}\cdots v_{j}^{n_j}I}}\times \mathcal{M}_{A/v_1^{n_1}\cdots v_{j}^{n_j}I}
		 \\& \xrightarrow{\mathrm{id} \times \mathrm{lev^*}} E_{A/v_1^{n_1}\cdots v_{j}^{n_j}I} \underset{\mathcal{M}_{A/v_1^{n_1}\cdots v_{j}^{n_j}I}}\times E_{A/v_1^{n_1}\cdots v_{j}^{n_j}I}
		\xrightarrow{\text{ group law }} E_{A/v_1^{n_1}\cdots v_{j}^{n_j}I}.
\end{align*}
And take care that
\[
	\phi_{n_j}'(\pi_{v_j}^{-n_j}) = \phi_{n_j}'(\mathrm{Frob}_{v_j}^{-n_j}(1)) = \mathrm{Frob}_{v_j}^{-n_j}\phi_{n_j}'(1).
\]
On the other hand what will $\phi_{n_j}(0 , \pi_{v_j}^{-n_j})^* m_{n_j}^*$ be? It is easy to show that 
\begin{align*}\begin{array}{cccccc}
	m_{n_j} \phi_{n_j}(0 , \pi_{v_j}^{-n_j}) :	&\mathcal{M}_{A/v_1^{n_1}\cdots v_{j}^{n_j}I}	&\rightarrow	&m_{n_j}^* E_{A/v_1^{n_1}\cdots v_{j -1}^{n_j -1}I}	&\rightarrow	&E_{A/v_1^{n_1}\cdots v_{j -1}^{n_j -1}I}		\\
							&\rotatebox{90}{$\cong$}	&		&	&			&\rotatebox{90}{$\cong$}	\\[-4pt]
							&\mathrm{Spec}\,\mathcal{O}_{w_i^j}(n) 	& &	&			&\mathrm{Spec}\,\mathcal{O}_{w_i^j}[T]
\end{array}\end{align*}
so we have
\begin{align*}\begin{array}{cccc}
	\phi_{n_j}(0 , \pi_{v_j}^{-n_j})^*m_{n_j}^* :	&\mathcal{O}_{w_i^j}[T]		&\rightarrow	&\mathcal{O}_{w_i^j}(n)		\\
								&\rotatebox{90}{$\in$}	&			&\rotatebox{90}{$\in$}		\\[-4pt]
								&T 					&\mapsto	&\xi_{n_j}
\end{array}\end{align*}
This asserts that it is the map which substitutes the value of the point for $T$. 
We know that there exists the Coleman power series for $u = (\sie{A/v_1^{n_1}\cdots v_{j}^{n_j}I})_{n_j} \in \varprojlim_{n_j} \mathcal{O}_{w_i^j}$. 
By summarizing, we have
\[
	\phi_{n_j}(0 , \pi_{v_j}^{-n_j})^* m_{n_j}^* t_{\phi_{n_j}'(\pi_{v_j}^{-n_j})} \vartheta_{E_{A/v_1^{n_1}\cdots v_{j}^{n_j}I}/\mathcal{M}_{A/v_1^{n_1}\cdots v_{j}^{n_j}I}}
	= \sie{A/v_1^{n_1}\cdots v_{j}^{n_j}I , \{ 1 \}} = \mathrm{Col}_{u} (\xi_{n_j}).
\]
Especially when $n_j = 0$,
\[
	\sie{A/v_1^{n_1}\cdots v_{j -1}^{n_j -1}I , \{ 1 \}} = \mathrm{Col}_{u}(0).
\]
Thus we obtain the following proposition:

\begin{Prop}
	\[
	t_{\phi_{n_j}'(\pi_{v_j}^{-n_j})} \vartheta_{E_{A/v_1^{n_1}\cdots v_{j}^{n_j}I}/\mathcal{M}_{A/v_1^{n_1}\cdots v_{j}^{n_j}I}} = \mathrm{Col}_{u}(T)
	\]
	Especially the constant terms of the Coleman power series is equal to the Siegel unit
	$\sie{A/v_1^{n_1}\cdots v_{j -1}^{n_j -1}I , \{ 1 \} } = \sie{(\bfvec{v}_{j-1} , \psi_{\bfvec{v}_{j-1}})}$. 
\end{Prop}

At the end of this section, we start to prove Theorem \ref{fsc}. \\
Let $(\bfvec{v}_r,\psi_{\bfvec{v}_r}) = \{(v_1,\psi_{v_1}),\ldots,(v_r,\psi_{v_r})\}$ be an element in $\Psi_M$ and 
take an element $(v,\psi_v) \in \Psi_M$. If $v \neq v_i$ for any $i=1,\ldots,r$, it is clear that
\[
	[\kappa_{(\bfvec{v}_r,\psi_{\bfvec{v}_r})}]_v = 0
\]
since $v_i$ is unramified in $\fd{v_1^{n_1}\cdots v_r^{n_r}}/\fd{}$. We will show the other case.
We may assume without loss of generality that $v = v_r$.
Let $w_r$ denote a prime of $\fd{}$ above $v_r$, $\fd{}_{w_r}$ the completion of $\fd{}$ at $w_r$, $w_r^{(i)}$ the prime of 
$E(\bfvec{v}_{i},\psi_{\bfvec{v}_{i}})$ above $w_r$ and $E_{w_r^{(i)}}$ the completion of $E(\bfvec{v}_i,\psi_{\bfvec{v}_i})$ 
at $w_r^{(i)}$ for any $i = 1,\ldots,r$. Put
\begin{align*}
	u_{E_{w_r^{(r)}}} &:= \displaystyle\prod_{i=1}^{r-1}\bfvec{D}_i \sie{(\bfvec{v}_r,\psi_{\bfvec{v}_r})}, \\
	u_{E_{w_r^{(r-i)}}} &:= N_{E_{w_r^{(r)}}/E_{w_r^{(r-i)}}}(u_{E_{w_r^{(r)}}}).
\end{align*}
In the same way of proof of Lemma \ref{3lem}, we have $u_{E_{w_r^{(0)}}} = 1$.
Moreover we have
\begin{align*}
	u_{E_{w_r^{(r-1)}}} &= N_{E_{w_r^{(r)}}/E_{w_r^{(r-1)}}}(u_{E_{w_r^{(r)}}})\\
		&= \bfvec{N}_r\displaystyle\prod_{i=1}^{r-1}\bfvec{D}_i \sie{(\bfvec{v}_r,\psi_{\bfvec{v}_r})}\\
		&= (1-\mathrm{Frob}_{v_r})\displaystyle\prod_{i=1}^{r-1}\bfvec{D}_i \sie{(\bfvec{v}_{r-1},\psi_{\bfvec{v}_{r-1}})}\\
		&\equiv 0 \quad \mathrm{mod}\,(E_{w_r^{(r-1)}})^M.
\end{align*}
So there exists an element $a \in E_{w_r^{(r-1)}}$ such that $a^M = u_{E_{w_r^{(r-1)}}}$.
Let $\chi$ be a character $G \rightarrow \mathbb{Q}/\mathbb{Z}$ such that the following diagram is commutative:
	\begin{align*}
	\begin{xy}
		(-1, 0)*{\psi_{v_r} :}, 
		(0, 25)*{\chi \ :}, 
		(10, 10)*{\fd{}_w^\times}, 
		(8, 0)*{\mathcal{O}_{v_r}^\times},
		(5, 5)*{\rotatebox{90}{$\cong$}},
		(5, 25)*{G_r}, 
		(40, 0)*{\mathbb{Z}/M\mathbb{Z}}, 
		(40, 25)*{\mathbb{Q}/\mathbb{Z}}, 
		\ar @{->} (12,0);(34,0), 
		\ar @{->}^{\mathrm{rec}} (5,13);(5,22), 
		\ar @{->}_{M^{-1}} (40,3);(40,22), 
		\ar @{->} (8,25);(34,25), 
	\end{xy}
	\end{align*}
Applying Corollary \ref{16} to $u = (u_{E_{w_r^{(n)}}})_n$, we have
\[
	(u_{E_{w_r^{(r)}}}/a , \chi)_{E_{w_r^{(r)}}/E_{w_r^{(r-1)}}} = \chi \circ \mathrm{rec}(\mathrm{Col}_u(0)/b^M)
\]
where $b$ is the element in $E_{w_r^{(r-1)}}^\times$ such that $a = \mathrm{Frob}_{v_r}(b) / b$.
By definition, the element $\kappa_{(\bfvec{v}_r,\psi_{\bfvec{v}_r})}$ is of the form
\[
	\kappa_{(\bfvec{v}_r,\psi_{\bfvec{v}_r})} = \displaystyle\prod_{i=1}^r\bfvec{D}_i\sie{(\bfvec{v}_r ,  \psi_{\bfvec{v}_r})}\cdot c
\]
where $c$ is an element in $(E_{w_r^{(r-1)}}^\times)^M$ which holds 
\[
	(1-\sigma_r)c = \displaystyle\frac{\displaystyle\prod_{i=1}^{r-1}\bfvec{D}_i \sie{(\bfvec{v}_r,\psi_{\bfvec{v}_r})}}
	{(1-\mathrm{Frob}_{v_r})\displaystyle\prod_{i=1}^{r-1}\bfvec{D}_i \sie{(\bfvec{v}_{r-1},\psi_{\bfvec{v}_{r-1}})}}
	= \displaystyle\frac{u_{E_{w_r^{(r)}}}}{a}.
\]
Remark that $\displaystyle\prod_{i=1}^r\bfvec{D}_i\sie{(\bfvec{v}_r ,  \psi_{\bfvec{v}_r})}$ is a unit, it is clear that
\[
	[\kappa_{(\bfvec{v}_r,\psi_{\bfvec{v}_r})}]_{v_r} = [c]_{v_r} = \underset{\lambda|v_r}\bigoplus\mathrm{val}_{v_r}(c)\lambda
\]
and we have $(u_{E_{w_r^{(r)}}}/a , \chi)_{E_{w_r^{(r)}}/E_{w_r^{(r-1)}}} = \mathrm{val}_{v_r}(c)$ since the definition.
The above implies that $(u_{E_{w_r^{(r)}}}/a , \chi)_{E_{w_r^{(r)}}/E_{w_r^{(r-1)}}}$ is equal to the coefficient of the $w_r$-part of 
$[\kappa_{(\bfvec{v}_r,\psi_{\bfvec{v}_r})}]_{v_r}$.
On the other hand, using Corollary \ref{20}, we have
\[
	\mathrm{Col}_u(0)/b^M = \kappa_{(\bfvec{v}_{r-1},\psi_{\bfvec{v}_{r-1}})}
\]
and then via the definition of $\chi$, we have
\[
	\chi\circ\mathrm{rec}(\mathrm{Col}_u(0)/b^M) = \psi_{v_{r}}(\kappa_{(\bfvec{v}_{r-1},\psi_{\bfvec{v}_{r-1}})}).
\]
This completes the proof of Theorem \ref{fsc}.

\subsection{The existence of prime}
The above, we consider a pair $(v,\psi_v) \in \Psi_M$. We now prove that there exists a ``good'' pair in our case. 
Let $\mathfrak{C}$ be a $p$-part of ideal class group of $\fd{}$.

\begin{Thm}\label{cheb}
	Given $\mathfrak{c} \in \mathfrak{C}$, a finite free $\mathbb{Z}/M\mathbb{Z}$-submodule $W$ of  $\fd{}^{\times} / (\fd{}^\times)^M$, 
	and surjection
	\[
	\varphi : W \rightarrow (\mathbb{Z}/M\mathbb{Z})[\Gal{\fd{}}{F}]
	\]
	compatible with the actions of the Galois group. Then there exists a prime $w$ of $\fd{}$ such that
	\begin{enumerate}
		\item $w \in \mathfrak{c}$
		\item there exist a prime $v$ of $F$ below $w$ and a map
			$\psi_v : \mathcal{O}_{F_v}^\times \rightarrow \mathbb{Z}/M\mathbb{Z}$ such that
		\begin{description}
			\item[(i)] $v$ splits completely in $\fd{}/F$,
			\item[(ii)] $[y]_v = 0$ for any $y \in W$,
			\item[(iii)] $\psi_v(w) = \varphi(w). $ 
		\end{description}
	\end{enumerate}
\end{Thm}

\begin{Pf}
	Let $H$ be the maximal unramified abelian $p$-extension of $\fd{}$, and regard $\mathfrak{C}$ as $\Gal{H}{\fd{}}$ by the global class field theory.
	\begin{align*}\begin{array}{ccc}
		\mathfrak{C}					&\xrightarrow{\sim}	&\Gal{H}{\fd{}}				\\
		\rotatebox{90}{$\in$}	&					&\rotatebox{90}{$\in$}		\\[-4pt]
		(\text{class of } \lambda)	&\mapsto			&\mathrm{Frob}_\lambda
	\end{array}\end{align*}
	where $\lambda$ is a prime of $\fd{}$. Let us consider a prime $w$ in $\mathfrak{c}$. It is clear that there exists a prime $v$ below $w$ such that 
	it holds (2)-(i). Let us consider this prime $v$. Since $W$ is a finite set, there exist only finitely many 
	primes which divide $(y) (y \in W)$, therefore we can take $v$ different from this and not dividing $(y)$.
	Now we construct $\psi_v$ satisfying (iii) from $\varphi$.
	First we show that $W$ is embedded to $\mathcal{O}_w^\times / (\mathcal{O}_w^\times)^M$ by an injection.
	For any $y \in W$, $y'$ denotes the image under the lift up to $\fd{}^\times$ in a completion of $\fd{},_w$.  
	For any $y \in W \backslash W^p$, if there exists an element $y'$ which is not in $(\mathcal{O}_w^\times)^p$, $W$ is 
	embedded by an injection in $\mathcal{O}_w^\times / (\mathcal{O}_w^\times)^M$ due to the claim(ii).
	But the condition needs
	\[
		y^{p^d} - y \not\equiv 0 \ \mathrm{mod}\ v^2.
	\]
	Here we use the fact that $y^{p^d}-y$ can be divided by $v$ at least 1 time where $p^d$ is the order of residue field of
	$F_v$. This requires that the morphism $\mathfrak{C} \rightarrow \mathbb{P}_{\mathbb{F}_q}^1$ induced by $y'$ is unramified at $v$.
	Since the assumption that $W$ is free $\mathbb{Z}/M\mathbb{Z}$-module and  $y' \in W \backslash W^p$,
	it is show that $y' \not\in F^p_ v$. Therefore $W$ is embedded in $\mathcal{O}_w^\times / (\mathcal{O}_w^\times)^M$ 
	under an injection with the finite exception of $v$. 
	Let $w,v$ be as above (and replace if we need). Fix an isomorphism $\mathcal{O}_w \cong \mathcal{O}_{F_v}$.
	There exists a homomorphism $\varphi'$ such that the following diagram holds: 
	\begin{align*}
	\begin{xy}
		(0, 0)*{W}, 
		(45, 0)*{\mathcal{O}_w^\times/(\mathcal{O}_w^\times)^M}, 
		(57, 0)*{\cong}, 
		(70, 0)*{\mathcal{O}_{F_v}^\times/(\mathcal{O}_{F_v}^\times)^M}, 
		(50, -10)*{\circlearrowright}, 
		(63, -20)*{(\mathbb{Z}/M\mathbb{Z})[\Gal{\fd{}}{F}]}, 
		\ar @{^{(}->}^{\text{injection}} (5,0);(35,0), 
		\ar @{->>}_{\varphi} (3,-3);(41,-17), 
		\ar @{.>}^{^{\exists}\varphi'} (70,-3);(70,-17), 
	\end{xy}
	\end{align*}
	Since $(\mathbb{Z}/M\mathbb{Z})[\Gal{\fd{}}{F}]$ is injective as $\mathbb{Z}/M\mathbb{Z}$-module.
	We regard naturally $\varphi'$ as $\mathcal{O}_{F_v}^\times \rightarrow (\mathbb{Z}/M\mathbb{Z})[\Gal{\fd{}}{F}]$
	with $\varphi'((\mathcal{O}_{F_v}^\times)^M) = 0$. Let us write it also $\varphi'$. 
	Let $\varrho : (\mathbb{Z}/M\mathbb{Z})[\Gal{\fd{}}{F}] \twoheadrightarrow \mathbb{Z}/M\mathbb{Z}$
	be a surjective $\mathbb{Z}/M\mathbb{Z}$-homomorphism such that 
	\begin{align*}\begin{array}{ccc}
		(\mathbb{Z}/M\mathbb{Z})[\Gal{\fd{}}{F}]		&\rightarrow				&\mathbb{Z}/M\mathbb{Z}		\\
		\rotatebox{90}{$\in$}						&						&\rotatebox{90}{$\in$}			\\[-4pt]
		\text{ unit of }\Gal{\fd{}}{F}					&\mapsto				&1							\\
		\text{ not unit of }\Gal{\fd{}}{F}				&\mapsto				&0
	\end{array}\end{align*}
	Put $\psi_v := \varrho\circ\varphi'$. This is the homomorphism what we need. 
\end{Pf}


\section{Proof of Iwasawa Main Conjecture}\label{lastsection}
In this section we describe the $p$-part of the ideal class group in function field (Theorem \ref{imc}) using Theorem \ref{fsc} 
and \ref{cheb} under some technical assumptions. We will prove this by imitating the proof of Iwasawa main conjecture in algebraic 
number fields. 

\subsection{Preparation of Iwasawa theory and Iwasawa main conjecture}\label{lastnotation}
Let us introduce some notations that we will use in Iwasawa theory. Let $F, p, p_\infty, A$ be as in section \ref{start}.
Let $I$ be an ideal of $A$. Put $K:=\fd{}$ for simplicity. Let $A_K$ be the integral closure of $A$ in $K$. Fix a prime $q$ of $A$
not dividing $I$. (Until the previous section we denoted by $q$ the order of the residue field of  $F$. Here we have changed the notation.)
Put $I_n := q^{n+1}I , K_n:=\fd{q^{n+1}}\ (n\geq0)$ for simplicity. Let
$K_\infty$ denotes $\underset{n\geq0}{\bigcup}\, K_n$, $\Kinfty$ a subfield of $K_\infty$ whose the Galois group
$\Gamma := \Gal{\Kinfty}{K_0}$ is isomorphic to $\mathbb{Z}_p$ and $K_n^{(p)} := (\Kinfty)^{\Gamma_n} , \Gamma_n := \Gamma^{p^n}$.
Put $\Delta := \Gal{K_0}{F}$. Assume that the order of $\Delta$ is prime to $p$.
Let $(K_n)_{q'}$ denote a completion of $K_n$ at the prime $q'$ of $K_n$ above $q$. Define symbols in Iwasawa theory as follows:
\begin{description}
	\item[$C_n$] 						: The $p$-part of the ideal class group of $K_n$.
	\item[$E_n$] 						: The unit group of $K_n$.
	\item[$\mathcal{E}_n \subset E_n$]	: The group generated by Siegel units of $K_n$.
	\item[$U_n$] 						: The unit group of $\prod(K_n^{(p)})_{q'}$, the product is over all the primes $q'$ of $(K_n^{(p)})$ above $q$.
	\item[$\overline{E_n}$] 				: The closure of $E_n \cap U_n$ in $U_n$.
	\item[$\overline{\mathcal{E}_n}$] 	: The closure of $\mathcal{E}_n \cap U_n$ in $U_n$.
\end{description}
And put $X_\infty := \varprojlim X_n$ for $X = C, E, \mathcal{E}, U, \overline{E}, \overline{\mathcal{E}}$,
where the inverse limit is taken with respect to the norm map. Especially $C_\infty$ is equal to the $p$-part of the ideal class
group of $\Kinfty$. Let $\chi$ denote a $p$-adic character in $\Delta$.
(Until the previous section we denoted by $\chi$ any character in $\Gal(H_\infty)(H)$. Here we have changed the notation.) 
For a $p$-adic character $\chi$, put
\[
e(\chi) := \frac{1}{\sharp\Delta}\displaystyle\sum_{\delta\in\Delta}\, \chi (\delta)^{-1}\delta.
\]
Define the Iwasawa algebra $\Lambda$ as 
\[
\Lambda := \mathbb{Z}_p[{\rm Im}\, \chi][[\Gamma]]. 
\]
Let $e(\chi)C_n ,  e(\chi)E_n ,  e(\chi)U_n ,  e(\chi)\overline{\mathcal{E}_n}$ denote 
$C_n(\chi) ,  E_n(\chi) ,  U_n(\chi) , \overline{\mathcal{E}_n}(\chi)$ for short. Similarly when $n = \infty$.
Then the following proposition is proved in \cite{OS}.

\begin{Prop}
	$C_\infty(\chi) ,  E_\infty(\chi) / \overline{\mathcal{E}_\infty}(\chi)$ are finitely generated $\Lambda$-modules.
\end{Prop}

Now we impose the following assumption on $C_\infty(\chi) , E_\infty(\chi)$ and $\overline{\mathcal{E}_\infty}(\chi)$.

\begin{Hyp}\label{hyp0}
	$C_\infty(\chi) ,  E_\infty(\chi)/\overline{\mathcal{E}_\infty}(\chi)$ are torsion $\Lambda$-modules.
	This means that $\Gamma$ satisfies the following conditions:\\
	Let $m_\Lambda$ be the maximal ideal of $\Lambda$. There exist elements $f_C,f_E \in m_\Lambda$ such that 
		\begin{align*}
		\sharp(C_\infty(\chi)/f_C (C_\infty(\chi)) < \infty \quad\text{ and } \\
		\sharp((E_\infty(\chi)/\overline{\mathcal{E}_\infty}(\chi))/f_E(E_\infty(\chi)/\overline{\mathcal{E}_\infty}(\chi))) < \infty.
	 \end{align*}
\end{Hyp}

From now on we assume hypothesis \ref{hyp0}. Now we recall the notion of pseudo-isomorphism of finitely generated torsion 
$\Lambda$-modules and that of characteristic ideals.

\begin{Def}[pseudo-isomorphic]
	Let $N,N'$ be finitely generated torsion $\Lambda$-modules. We say that $N$ and $N'$ are {\bf pseudo-isomorphic} when 
	there exists a homomorphism $N \rightarrow N'$ of $\Lambda$-modules whose kernel and cokernel are finite.
	If $N$ and $N'$ are pseudo-isomorphic then we write $N \sim N'$.
 \end{Def}

The following remark is a well-known fact about finitely generated torsion $\Lambda$-modules.

\begin{Rmk}
	For any finite generated torsion $\Lambda$-module $N$, there are finitely many
	$f_i\ (i=1, \ldots, r) \in m_\Lambda$ unique up to rearrangement such that
	\[
	N \sim \displaystyle\bigoplus_{i=1}^r\, \Lambda/f_i\Lambda.
	\] 
\end{Rmk}

\begin{Def}[characteristic ideal]
	For the above $f_i$, put $f := \displaystyle\prod_{i} f_i$. 
	We call the ideal $f\Lambda$ of $\Lambda$ the
	{\bf characteristic ideal} of $N$, and write ${\rm char}\, (N)$.
\end{Def}

In this setting, we have the following.

\begin{Thm}[Iwasawa main conjecture]\label{imc}
	For any $p$-adic irreducible character $\chi$ in $\Delta$,
	\begin{align*}
		{\rm char}\, (C_\infty(\chi)) = {\rm char}\left(E_\infty(\chi)/\overline{\mathcal{E}_\infty}(\chi)\right).
	 \end{align*}
\end{Thm}

The goal of this section is to prove this theorem. We will prove it in section \ref{last} 

\subsection{Iwasawa theory}

First we write $a_n \approx b_n$ when $a_n/b_n$ is bounded from below and above independently of $n$, where $a_n$ and $b_n$
are sequences of positive integers.

\begin{Prop}\label{hyp1}
	For any $n\geq0$, $\sharp(C_n) \approx [\overline{E_n}:\overline{\mathcal{E}_n}]$. 
 \end{Prop}
\begin{Pf}
We refer the reader to \cite[Appendix-Lemma 6.6]{La}. The equality between
the index number of unit groups and ideal class number in global function fields is proved by L. Yin \cite{Yi}.

\end{Pf}

Next for any $n\geq0$, put $\Gamma_n := \Gal{\Kinfty}{K_n}$. 
Let $J_n$ be an ideal of $\Lambda$ generated by $\{\gamma-1\mid\gamma\in\Gamma_n\}$.
If $\gamma$ is a generator of $\Gamma$, we have $J_n = (\gamma^{p^n}-1)\Lambda$.
Put $\Lambda_n := \Lambda/J_n\Lambda$. Then we have $\Lambda_n \cong \mathbb{Z}_p[\mathrm{Im} \chi][\Gal{K_n}{K_0}]$.
For $\Lambda$-module $N$, let $N_{\Lambda_n}$ denote a tensor product $N \displaystyle\otimes_{\Lambda} \Lambda_n$
and regard it as a $\Lambda_n$-module.
Let us consider the pseudo-isomorphisms described above: 
\begin{eqnarray*}
	(C_\infty(\chi))&\sim&\displaystyle\bigoplus_{i=1}^k\, \Lambda/f_i\Lambda , \\
	\left(E_\infty(\chi)/(\overline{\mathcal{E}_\infty}(\chi))\right)&\sim&\displaystyle\bigoplus_{i=1}^l\, 
	\Lambda/h_i\Lambda.
\end{eqnarray*}
For each $f_i , h_i$ corresponding to the above, we define
\begin{eqnarray*}
	f_\chi &:=& \displaystyle\prod_{i=1}^k\, f_i , \\
	h_\chi &:=& \displaystyle\prod_{i=1}^l\, h_i.
\end{eqnarray*}
It is clear that ${\rm char}(C_\infty(\chi)) = f_\chi\Lambda$ and ${\rm char}\left(E_\infty(\chi)/(\overline{\mathcal{E}_\infty}(\chi))\right) = h_\chi\Lambda$.

\begin{Lem}\label{lem1}
	Assume that $a_1,  a_2 \in \Lambda ,  a_1|a_2 , 
	\sharp((\Lambda/a_1\Lambda)_{\Lambda_n}) \approx \sharp((\Lambda/a_2\Lambda)_{\Lambda_n})$. Then
	\[ a_1\Lambda = a_2\Lambda.  \]
\end{Lem}
\begin{Pf}
	We refer the reader to \cite[Appendix-Corollary7.3]{La}.
\end{Pf}

\begin{Lem}\label{lem2}
	Let $\chi$ be a non trivial character of $\Delta$. Then there exists an ideal $\mathscr{A}$ of  $\Lambda$ which is
	 of finite index such that followings hold.
	\begin{enumerate}
		\item For any $n\geq0$, there exist elements $\mathfrak{c}_1 ,  \ldots ,  \mathfrak{c}_k \in C_n(\chi)$ such that
				\[ \mathscr{J}(D_q)\mathscr{A}\mathrm{Ann}(\mathfrak{c}_i) \subseteq f_i\Lambda_n. \]
		\item For any $\eta \in \mathscr{J}(D_q)\mathscr{A}$ and $n\geq0$, there exists a $\Lambda$ homomorphism
			$\theta_{n, \eta}:\overline{E_n}(\chi)\rightarrow\Lambda_n$ such that
				\[ \eta^4 h_\chi\Lambda_n \subseteq \theta_{n, \eta}(\overline{\mathcal{E}_n}(\chi)). \]
	\end{enumerate}
	Where $D_q$ is the decomposition group of $q$ in $\Gal{K_\infty^{(p)}}{F}$ and $\mathscr{J}(D_q)$
	the ideal of the Iwasawa algebra of $\Gal{K_\infty^{(p)}}{F}$ generated by the elements $\sigma -1, \sigma \in D_q$.
\end{Lem}

\begin{Pf}
	We refer the reader to \cite[Proposition 6.5, Corollary 7.10]{R2}.
\end{Pf}

\subsection{Proof of Theorem \ref{imc}}\label{last}
In this section we fix a positive integer $n$.
Let $\lambda'$ be a prime of $F$ which splits completely in $K_n/F$. Recall that for $y \in K_n^\times$, 
we denote $(y)_{\lambda'}\in\mathscr{I}_{\lambda'}$ the $\lambda'$-part of $(y)\in\mathscr{I}$, 
$[y]_{\lambda'}\in\mathscr{I}_{\lambda'}/M\mathscr{I}_{\lambda'}$ the image under the projection where $M$ is a power of $p$.
Let $\lambda$ be a prime of $K_n$ above $\lambda'$.
Then $\mathscr{I}_{\lambda'}(\chi) := e(\chi)\mathscr{I}_{\lambda'}$ is a free $\Lambda_n$ module of degree one over
$\Lambda_n \cong \mathbb{Z}_p[\mathrm{Im}\,\chi][\Gal{K_n}{K_0}]$, generated by 
$\lambda(\chi) := e(\chi)\lambda$.
We define the map $\nu_{\lambda ,  \chi} : K_n^\times \rightarrow \Lambda_n \cong \mathscr{I}_{\lambda'}(\chi)$ satisfying
\[
	\nu_{\lambda ,  \chi}(y)\lambda(\chi) = e(\chi)(y)_{\lambda'}
\]
and the map $\overline{\nu_{\lambda ,  \chi}} : K_n^\times/(K_n^\times)^M \rightarrow \Lambda_n/M\Lambda_n$ satisfying
\[
	\overline{\nu_{\lambda ,  \chi}}(y)\lambda(\chi) = e(\chi)[y]_{\lambda'}.
\]
From the section \ref{kap}, for the set $\Psi_M$, $(\bfvec{v} ,  \psi_{\bfvec{v}}) \in \Psi_M$ and the field extension $K_n/F$, 
there exist Kolyvagin's derivative classes $\kappa_{(\bfvec{v} ,  \psi_{\bfvec{v}})} \in K_n^\times/(K_n^\times)^M$.
Be careful that the ideal $I$ in the section \ref{kol} is replaced with $q^{n+1}I$ here, similarly $\fd{}$ replaced with $K_n$. 
We prepare two lemmas.

\begin{Lem}\label{lem4}
	Let $M$ be a sufficiently large power of $p$ and $(\bfvec{v} ,  \psi_{\bfvec{v}})\in\Psi_M$. 
	Take a prime $\lambda'$ of $F$ which divides $\bfvec{v}$ and a prime $\lambda$ of $K_n$ above $\lambda'$.
	Let $B_n$ be a subgroup of $C_n$ generated by $\bfvec{v}/\lambda'$.
	Let $\mathfrak{c}$ denote the class of ideal $e(\chi)\lambda$, $W$ a sub $\Lambda_n$-module of 
	$K_n^\times/(K_n^\times)^M$ generated by $e(\chi)\kappa_{(\bfvec{v} ,  \psi_{\bfvec{v}})}$.
	Assume that $\eta ,  a\in\Lambda_n$ satisfy the following conditions:
	\begin{enumerate}
		\item $\Lambda_n/a\Lambda_n$ is finite.
		\item ${\rm Ann}(\mathfrak{c}) \subseteq \Lambda_n$ satisfies
				$\eta{\rm Ann}(\mathfrak{c}) \subseteq a\Lambda_n$ in $C_n(\chi)/B_n(\chi)$.					
	\end{enumerate}
	Then there exists a homomorphism 
	$\varphi : W \rightarrow \Lambda_n/M\Lambda_n \cong (\mathbb{Z}/M\mathbb{Z})[\Gal{K}{F}]$ 
	which keeps the action of Galois group such that
	\[ a\varphi(e(\chi)\kappa_{(\bfvec{v} ,  \psi_{\bfvec{v}})}) =
				 \eta\overline{\nu_{\lambda, \chi}}(\kappa_{(\bfvec{v} ,  \psi_{\bfvec{v}})}). \]
\end{Lem}

\begin{Pf}
	Let $\beta$ be one of a lifting of $e(\chi)\kappa_{(\bfvec{v} ,  \psi_{\bfvec{v}})}$ to $K_n^\times$.
	Via the definition of $\nu_{\lambda, \chi}$, we have
	\begin{eqnarray*}
	e(\chi)(\beta) &=& e(\chi)(\beta)_{\lambda'} + \displaystyle\sum_{\lambda''\neq\lambda'}e(\chi)(\beta)_{\lambda''} \\
			&=& \nu_{\lambda ,  \chi}(\beta)\lambda(\chi) + \displaystyle\sum_{\lambda''\neq\lambda'}e(\chi)(\beta)_{\lambda''}.
	\end{eqnarray*}
	Here Theorem \ref{fsc} implies that when $\lambda''$ does not divide $\bfvec{v}$, 
	$(\beta)_{\lambda''} \in M\mathscr{I}_{\lambda''}$. Since we assume that $M$ is a sufficiently large power of $p$, $M$
	annihilates $C_n(\chi)$.
	Therefore $\nu_{\lambda ,  \chi}(\beta)\lambda(\chi) \in {\rm Ann}(\mathfrak{c})$ since 
	$\nu_{\lambda ,  \chi}(\beta)\lambda(\chi)$ is zero in $C_n(\chi)/B_n(\chi)$ by the definition of $B_n(\chi)$.
	The second assumption implies
	\[
	\nu_{\lambda ,  \chi}(\beta)\lambda(\chi) \in a\Lambda_n. 
	\]
	The first assumption implies that $\nu_{\lambda ,  \chi}(\beta)\lambda(\chi) /a$ is well-defined. Let us write it by $\sigma$.
	Recall that $W$ is generated by  $e(\chi)\kappa_{(\bfvec{v} ,  \psi_{\bfvec{v}})}$ over 
	$\Lambda_n \cong \mathbb{Z}_p[\Gal{K_n}{K_0}]$. Define a homomorphism 
	$\varphi : W \rightarrow \Lambda_n/M\Lambda_n$ by
	\[
	\varphi\left(y(e(\chi)\kappa_{(\bfvec{v} ,  \psi_{\bfvec{v}})})\right) = y\sigma 
	\]
	for any $y \in \Lambda_n$. It is clear that this keeps the Galois action. Therefore we check that it is well-defined
	and is independent of the way of taking $\beta$. 
	Assume that $y e(\chi)\kappa_{(\bfvec{v} ,  \psi_{\bfvec{v}})} = 0$ i.e., there exists an element $x$ in $\Lambda_n$ such that
	$y\beta = x^M$. Especially $y[e(\chi)\kappa_{(\bfvec{v} ,  \psi_{\bfvec{v}})}]_{\lambda'} = 0$.
	We have to prove $y\sigma \in M\Lambda_n$.
	We assume that $M$ is large enough, so we can assume $(M/\sharp C(\chi))(\mathscr{I}_{\lambda'}	
	(\chi)/M\mathscr{I}_{\lambda'}(\chi))$ is included in $\Lambda_n[e(\chi)\kappa_{(\bfvec{v} ,  \psi_{\bfvec{v}})}]$.
	Therefore we have $y \in (\sharp C_n(\chi))\Lambda_n$. Thus
	\begin{align*}
		e(\chi)(x) = \displaystyle\sum_{\lambda''}e(\chi)(x)_{\lambda''} & \\
				= M^{-1}e(\chi)(y\beta)_{\lambda'}
				 &+ \displaystyle\sum_{\lambda''\text{ divides }(\bfvec{v}/(\lambda'))}e(\chi)(x)_{\lambda''} \\
				 &+ \displaystyle\sum_{\lambda''\text{ does not divide }(\bfvec{v})}e(\chi)(y\beta)_{\lambda''} \\
				\equiv M^{-1}e(\chi)(y\beta)_{\lambda'} 
				 &+ \displaystyle\sum_{\lambda''\text{ divides }(\bfvec{v}/(\lambda'))}e(\chi)(x)_{\lambda''} \\
				&\quad {\rm mod}\ \sharp C_n(\chi)\mathscr{I}(\chi).
	\end{align*}
	Here $\sharp C_n(\chi)$ annihilates $C_n(\chi)$, hence $M^{-1}e(\chi)(y\beta)_{\lambda'}$ is equal to zero \\ in 
	$C_n(\chi)/B_n(\chi)$.
	Therefore $M^{-1}\nu_{\lambda, \chi}(y\beta)\mathfrak{c} = 0$ and we have
	\[
	y\sigma a = \eta\nu_{\lambda, \chi}(y\beta) \in Ma\Lambda_n.
	\]
	Dividing both sides by $a$, we have $y\sigma \in M\Lambda_n$.
\end{Pf}

\begin{Lem}\label{lem3}
	For any $p$-adic irreducible character $\chi$ of $\Delta$,  ${\rm char}(C_\infty(\chi))$ divides
	${\rm char}\left(E_\infty(\chi)/(\overline{\mathcal{E}_\infty}(\chi))\right)$, i.e.,
	\[ f_\chi | h_\chi.  \]
\end{Lem}

\begin{Pf}
	If $\chi = 1$, it is clear that 
	$C_n(\chi) = (\frac{1}{\sharp\Delta})\displaystyle\sum_{\delta\in\Delta}\chi(\delta)^{-1}\delta C_n = C_n$.
	Since $\gamma$ denotes the topological generator of $\Gamma$, it is easy to show that $C_n/C_n^{\gamma-1} = C_0$. 
	Since the definition and the fact that $C_0 = 1$, we have $C_n = C_n^{\gamma-1}$.
	Via the usual identification $\gamma-1 \mapsto X$, $\mathbb{Z}_p[[\Gamma]] \cong \mathbb{Z}_p[[X]]$ in Iwasawa theory,
	$C_n$ is regarded as a  $\mathbb{Z}_p[[X]]$-module. The above identification imply that
	\[ C_n = (X)C_n \]
	where $(X)$ is the maximal ideal of local ring, especially it is Jacobson's radical. By the Nakayama's lemma, 
	we have $C_n = 0$. Therefore $C_\infty = 0$ and $f_\chi$ is a unit.
	
	Assume that $\chi \neq 1$. Remark that $\overline{\mathcal{E}_n}(\chi)$ is generated by 
	$\sie{A/q^{n+1}I ,  \{1\}}(\chi) := e(\chi)\sie{A/q^{n+1}I ,  \{1\}}$.
	Let us consider an ideal $\mathscr{C}$ of $\Lambda$ which has finite index and elements 
	$\mathfrak{c}_1 ,  \ldots ,  \mathfrak{c}_k$  of $C_n(\chi)$, satisfying the conditions in Lemma \ref{lem2}.
	Moreover take another element $\mathfrak{c}_{k+1}$ in $C_n(\chi)$.
	(For example, $\mathfrak{c}_{k+1} = 0$). There exists an element $\eta \in \mathscr{C}$ such that $\Lambda_n/\eta\Lambda_n$ 
	is finite, (i.e. $\eta$ is prime to $\gamma^{p^n} -1$), and for there exists a homomorphism $\theta_{n, \eta}$ over $\Lambda$ such that 
	$\eta^4 h_\chi\Lambda_n \subseteq \theta_{n , \eta}(\overline{\mathcal{E}_n} (\chi))$.
	Here without loss of generality, we assume that
	\[
		\theta_{n , \eta}(\sie{A/q^{n+1}I , \{1\}}(\chi)) = \eta^4 h_\chi.
	\]
	Let $M\in \mathbb{Z}$ be a power of $p$ large enough. From now we use Theorem \ref{cheb} to construct 
	$\lambda_i$ satisfying the following conditions each $i$ with $1\leq i \leq k+1$, 
	\begin{enumerate}
		\item $\lambda_i \in \mathfrak{c}_i$,
		\item there exists a pairing $(\lambda'_i ,  \psi_{\lambda'_i}) \in \Psi_M$ where $\lambda'_i$ is a prime below to $\lambda_i$ which splits
			 completely in $K_n/F$ and $\psi_{\lambda'_i}$ is a continuous surjective homomorphism, 
		\item If $i = 1$ then $\overline{\nu_{\lambda_1 ,  \chi}}(\kappa_{(\lambda'_1 ,  \psi_{\lambda'_1})}) = \eta^4 h_\chi$
		\item If $2 \leq i \leq k+1$ then 
				$f_{i-1}\overline{\nu_{\lambda_i ,  \chi}}(\kappa_{(\bfvec{v}'_i ,  \psi_{\bfvec{v}'_i})})
				= \eta\overline{\nu_{\lambda_{i-1} ,  \chi}}(\kappa_{(\bfvec{v}'_{i-1} ,  \psi_{\bfvec{v}'_{i-1}})})$ \\
				where $(\bfvec{v}'_i ,  \psi_{\bfvec{v}'_i}) 
				= \{ (\lambda'_1 ,  \psi_{\lambda'_1}) ,  \ldots ,  (\lambda'_i ,  \psi_{\lambda'_i}) \}$. 
	\end{enumerate}
	First we will show to choose $\lambda_i$ by induction on $i$. When $i = 1$, 
	let $\iota : E_n(\chi)/E_n(\chi)^M \rightarrow \overline{E_n}(\chi)/\overline{E_n}(\chi)^M$ be a natural embedding.
	Applying Theorem \ref{cheb} to the case when $\mathfrak{c}=\mathfrak{c}_1 ,  W = (E/E^M)(\chi)$ and $\varphi$ as the composition of
	\[
		W \xrightarrow{\iota} \overline{E}(\chi)/\overline{E}(\chi)^M \xrightarrow{\theta_{n, \eta}} \Lambda_n/M\Lambda_n
		\xrightarrow{\cong} (\mathbb{Z}/M\mathbb{Z})[\Gal{K}{F}].
	\]
	Then there exist a prime $\lambda_1 \in \mathfrak{c}$ and a pairing $(\lambda'_1, \psi_{\lambda'_1})$
	where $\lambda'_1$ is below to $\lambda_1$.
	It is clear from the choice of $\lambda_1$ that (1) and (2) are satisfied.
	We will prove that (3) also holds. Theorem \ref{fsc} and \ref{cheb} imply
	\begin{align*}
		\overline{\nu_{\lambda_1 ,  \chi}}(\kappa_{(\lambda'_1 ,  \psi_{\lambda'_1})})\lambda_1(\chi)
			&= e(\chi)[\kappa_{(\lambda'_1 ,  \psi_{\lambda'_1})}]_{\lambda'_1} \\
			&= e(\chi)\psi_{\lambda'_1}(\kappa_\emptyset) \\
			&= \varphi(\kappa_{\emptyset})\lambda_1(\chi) \\
			&= \theta_{n, \eta}(\sie{\emptyset})\lambda_1(\chi) \\
			&= \eta h_\chi \lambda_1(\chi).
	\end{align*}
	Thus $\overline{\nu_{\lambda_1 ,  \chi}}(\kappa_{(\lambda'_1 ,  \psi_{\lambda'_1})}) = \eta h_\chi$.
	
	Continue this induction process, assume that we can choose $\lambda_1 ,  \ldots ,  \lambda_{i-1}$. Put 
	$(\bfvec{v}'_{i-1} ,  \psi_{\bfvec{v}'_{i-1}}) 
				= \{ (\lambda'_1 ,  \psi_{\lambda'_1}) ,  \ldots ,  (\lambda'_{i-1} ,  \psi_{\lambda'_{i-1}}) \}$.
	Let $W_i$ be a $\Lambda_n$-submodule in $K_n^\times/(K_n^\times)^M$ generated by 
	$e(\chi)\kappa_{(\bfvec{v}'_{i-1} ,  \psi_{\bfvec{v}'_{i-1}})}$. It follows from Lemma \ref{lem4} that
	there exists a map $\varphi_i : W_i \rightarrow \Lambda_n/M\Lambda_n$ such that 
	\[
		f_{i-1}\varphi_i(e(\chi)\kappa_{(\bfvec{v}'_{i-1} ,  \psi_{\bfvec{v}'_{i-1}})}) 
			= \eta\overline{\nu_{\lambda_{i-1} ,  \chi}}(\kappa_{(\bfvec{v}'_{i-1} ,  \psi_{\bfvec{v}'_{i-1}})}).
	\]
	Applying Theorem \ref{cheb} with $\mathfrak{c}=\mathfrak{c}_i ,  W=W_i ,  \varphi = e(\chi)\varphi_i$,
	there exist a prime $\lambda_i \in \mathfrak{c}$ and a pairing $(\lambda'_i ,  \psi_{\lambda'_i})$ where 
	$\lambda'_i$ is below to $\lambda_i$. Now we show that (4) holds. Similarly when $i =1$, we can calculate
	\begin{align*}
		f_{i-1}\overline{\nu_{\lambda_i ,  \chi}}(\kappa_{(\bfvec{v}'_i ,  \psi_{\bfvec{v}'_i})})
			&= f_{i-1}e(\chi)[\kappa_{(\bfvec{v}'_i ,  \psi_{\bfvec{v}'_i})}]_{\lambda_i} \\
			&= f_{i-1}\psi_{\lambda'_i}(e(\chi)\kappa_{(\bfvec{v}'_{i-1} ,  \psi_{\bfvec{v}'_{i-1}})}) \\
			&= f_{i-1}\varphi_i(e(\chi)\kappa_(\bfvec{v}'_{i-1} ,  \psi_{\bfvec{v}'_{i-1}}))\lambda_i(\chi)  \\
			&= \eta\overline{\nu_{\lambda_{i-1} ,  \chi}}(\kappa_{(\bfvec{v}'_{i-1} ,  \psi_{\bfvec{v}'_{i-1}})})\lambda_i(\chi).
	\end{align*}
	Therefore we have $f_{i-1}\overline{\nu_{\lambda_i ,  \chi}}(\kappa_{(\bfvec{v}'_i ,  \psi_{\bfvec{v}'_i})})
				= \eta\overline{\nu_{\lambda_{i-1} ,  \chi}}(\kappa_{(\bfvec{v}'_{i-1} ,  \psi_{\bfvec{v}'_{i-1}})})$.
	From the above we construct $\lambda_i$ inductively.
	
	What we have to prove is that $f_\chi$ divides $h_\chi$. It is easy to calculate that 
	\begin{align*}
		\eta^{k+4}h_\chi &= \eta^k \overline{\nu_{\lambda_1, \chi}}(\kappa_{(\bfvec{v}'_{1} ,  \psi_{\bfvec{v}'_{1}})}) \\
			&= \eta^{k-1}f_1\overline{\nu_{\lambda_2, \chi}}(\kappa_{(\bfvec{v}'_{2} ,  \psi_{\bfvec{v}'_{2}})}) \\
			&= \cdots \\
			&= \left( \displaystyle\prod_{j=1}^k f_j\right) 
				\overline{\nu_{\lambda_{k+1}, \chi}}(\kappa_{(\bfvec{v}'_{k+1} ,  \psi_{\bfvec{v}'_{k+1}})}),
	\end{align*}
	therefore $f_\chi = \displaystyle\prod_{j=1}^k f_j$ divides $\eta^{k+4}h_\chi$. Since our assumption of $\Gamma$,
	$char(C_\infty(\chi))$ is prime to $\mathscr{J}(D_q)$. Therefore $f_\chi$ divides $h_\chi$.
\end{Pf}

\begin{Pf}[of Theorem \ref{imc}]
	It is easy to show the theorem using Lemma \ref{lem1}, \ref{lem3}, and Proposition \ref{hyp1}. Put
	\begin{align*}
	&f := \displaystyle\prod_{\chi}f_\chi = \displaystyle\prod_{\chi}{\rm char}(C_\infty(\chi)) \\
	&h := \displaystyle\prod_{\chi}h_\chi = \displaystyle\prod_{\chi}
		{\rm char}\left(E_\infty(\chi)/(\overline{\mathcal{E}_\infty}(\chi))\right).
	\end{align*}
	For any $m\geq0$, we have
	\begin{align*}
		&\sharp(\Lambda_m/f\Lambda_m) \approx \displaystyle\prod_{\chi}\sharp(\Lambda_m/f_\chi\Lambda_m)
			\approx \displaystyle\prod_{\chi}\sharp(C_m(\chi)) = \sharp(C_m) \\
		&\sharp(\Lambda_m/h\Lambda_m) \approx \displaystyle\prod_{\chi}\sharp(\Lambda_m/h_\chi\Lambda_m)
			\approx \displaystyle\prod_{\chi}[\overline{E_m}(\chi) : \overline{\mathcal{E}_m}(\chi)]
			= [\overline{E_n} : \overline{\mathcal{E}_n}].
	\end{align*}
	Proposition \ref{hyp1} implies that $\sharp(\Lambda_m/f\Lambda_m) \approx \sharp(\Lambda_m/h\Lambda_m)$. 
	Lemma \ref{lem3} implies that $f|g$. Since $f,g \in \Lambda$, we have $f\Lambda = g\Lambda$ by using Lemma \ref{lem1}.
	Using Lemma \ref{lem3} once more, we have $f_\chi\Lambda = g_\chi\Lambda$ for any $\chi$. Therefore we have
	\[
	{\rm char}\, (C_\infty(\chi)) = {\rm char}\left(E_\infty(\chi)/(\overline{\mathcal{E}_\infty}(\chi))\right).
	\]
	This is what we need to prove.
\end{Pf}



\begin{thebibliography}{99}
	\bibitem[BL]{BL} F.\ Bars,\ I.\ Longhi.\quad {\em Coleman's power series and Wiles' reciprocity for rank 1 Drinfeld modules.}\quad Journal of Number Theory.\ {\bf 129}, 789--805 (2009).
	\bibitem[Co]{Co} R.\ F.\ Coleman.\quad {\em Division Values in Local Fields.}\quad Inventiones mathematicae\ {\bf 53}, 91--116 (1979).
	\bibitem[Dr]{Dr} V.\ G.\ Drinfeld.\quad {\em Elliptic modules.}\quad Math.\ USSR-Sbornik\ {\bf 23}, 561--592 (1973).
	\bibitem[Ha]{Ha} M.\ Hazewinkel.\quad {\em Corps de classes local. (appendix of Groupes alg$\acute{e}$briques.)}\quad Paris Masson (1970)
	\bibitem[K1]{K1} V.\ A.\ Kolyvagin.\quad {\em Euler system.}\quad Prog.\ Math.\ {\bf 87}, 435--483, Birkh$\ddot{\mathrm{a}}$user, Boston (1990).
	\bibitem[K2]{K2} V.\ A.\ Kolyvagin.\quad {\em The Mordell-Weil and Shafarevich-Tate groups for Weil elliptic curves.}\quad Math.\ USSR-Izvestiya\ {\bf 33}, 474--499 (1989).
	\bibitem[Ka]{Ka} K.\ Kato.\quad {\em p-adic Hodge theory and values of zeta functions of modular forms.}\quad Ast$\acute{\mathrm{e}}$risque {\bf 295}, 117--290 (2004).
	\bibitem[KY]{KY} S.\ Kondo, \ S.\ Yasuda.\quad {\em Zeta elements in the K-theory of Drinfeld modular varieties.}\quad Math.\ Ann. {\bf 354}, 529--587 (2012).
	\bibitem[La]{La} S.\ Lang.\quad {\em Cyclotomic Fields I and I\hspace{-.1em}I.}\quad Graduate Text in Mathematics, {\bf 121} Springer (1990).
	\bibitem[LLTT]{LLTT} K.\ F.\ Lai,\ I.\ Longhi,\ K.\ Tan,\ F.\ Trihan.\quad {\em The iwasawa main conjecture for semistable abelian varieties over function fields.}\quad Preprint, available at arXiv, \texttt{http://arxiv.org/abs/1406.6128}.
	\bibitem[OS]{OS} H.\ Oukhaba\ S.\ Vigui$\acute{{\rm e}}$\quad {\em The Gras Conjecture in function fields by euler systems.}\quad Bull.\ London Math.\ Soc.\ {\bf 43}, 523--535 (2011).
	\bibitem[R1]{R1} K.\ Rubin\quad {\em Euler systems.}\quad Annals of Mathematics studies {\bf 147}, Prinston Univ.\ Press (2000).
	\bibitem[R2]{R2} K.\ Rubin\quad {\em The ``main conjectures'' of Iwasawa theory for imaginary quadratic fields.}\quad Invent. Math. {\bf 103}, no.1, 25--68 (1991) 
	\bibitem[Yi]{Yi} L.\ Yin\quad {\em Index-class number formulas over global function fields.}\quad Compositio Math. {\bf 109}, 49--66 (1997).
\end{thebibliography}
\end{document}